\documentclass[12pt]{article}
\usepackage{amssymb}
\usepackage{latexsym}
\usepackage{graphics}
\newcommand{\ben}{\begin{enumerate}}
\newcommand{\een}{\end{enumerate}}
\newcommand{\ble}{\begin{lem}}
\newcommand{\ele}{\end{lem}}
\newcommand{\bth}{\begin{thm}}
\renewcommand{\eth}{\end{thm}}
\newcommand{\bpr}{\begin{prop}}
\newcommand{\epr}{\end{prop}}
\newcommand{\bco}{\begin{cor}}
\newcommand{\eco}{\end{cor}}
\newcommand{\bcon}{\begin{conj}}
\newcommand{\econ}{\end{conj}}
\newcommand{\bde}{\begin{defn}}
\newcommand{\ede}{\end{defn}}
\newcommand{\bex}{\begin{exa}}
\newcommand{\eex}{\end{exa}}
\newcommand{\barr}{\begin{array}}
\newcommand{\earr}{\end{array}}
\newcommand{\btab}{\begin{tabular}}
\newcommand{\etab}{\end{tabular}}
\newcommand{\beq}{\begin{equation}}
\newcommand{\eeq}{\end{equation}}
\newcommand{\bea}{\begin{eqnarray*}}
\newcommand{\eea}{\end{eqnarray*}}
\newcommand{\bce}{\begin{center}}
\newcommand{\ece}{\end{center}}
\newcommand{\bpi}{\begin{picture}}
\newcommand{\epi}{\end{picture}}
\newcommand{\bfi}{\begin{figure} \begin{center}}
\newcommand{\efi}{\end{center} \end{figure}}

\newcommand{\bsl}{\begin{slide}{}}
\newcommand{\esl}{\end{slide}}

\newcommand{\pf}{\noindent{\bf Proof}\hspace{7pt}}
\newcommand{\qed}{\rule{1ex}{1ex}}

\newcommand{\hso}[1]{\hspace{-1pt}}

\newcommand{\sbs}{\subset}
\newcommand{\sbe}{\subseteq}

\newcommand{\spe}{\supseteq}

\newcommand{\zh}{\hat{0}}

\def\<{\langle}
\def\>{\rangle}

\newcommand{\cA}{{\cal A}}

\newcommand{\cB}{{\cal B}}

\newcommand{\cC}{{\cal C}}

\newcommand{\cF}{{\cal F}}

\newcommand{\cH}{{\cal H}}

\newcommand{\cL}{{\cal L}}

\newcommand{\cO}{{\cal O}}

\newcommand{\cT}{{\cal T}}

\newcommand{\cX}{{\cal X}}
\newcommand{\cY}{{\cal Y}}
\newcommand{\cZ}{{\cal Z}}

\renewcommand{\bar}{\overline}

\newcommand{\inc}{\mathop{\rm inc}\nolimits}

\newcommand{\Rad}{\mathop{\rm Rad}\nolimits}

\newcommand{\rank}{\mathop{\rm rank}\nolimits}

\def\flexbox#1{\mathchoice{\mbox{#1}}{\mbox{#1}}{\mbox{\scriptsize #1}}%
{\mbox{\tiny #1}}}
\def\SL{\mathop{\flexbox{\rm SL}}\nolimits}
\def\PSL{\mathop{\flexbox{\rm PSL}}}

\newcommand{\SP}{\mathop{\flexbox{\rm Sp}}}

\newcommand{\Char}{\mathop{\flexbox{\rm Char}}}

\newcommand{\FF}{{\mathbb F}}

\newcommand{\NN}{{\mathbb N}}

\newcommand{\ZZ}{{\mathbb Z}}

\newcommand{\sform}{{\sf s}}

\renewcommand{\inc}{\star}
\newcommand{\typ}{\mathop{\rm typ}}%

\newcommand{\meet}{\wedge}

\newcommand{\dfn}{\em}

\newcommand{\after}{\mathbin{\circ}}
\newcommand{\GL}{\mathop{\rm GL}\nolimits}
\newcommand{\Sp}{\mathop{\rm Sp}\nolimits}

\newcommand{\Stab}{\mathop{\rm Stab}}

\newcommand{\AI}[1]{\item[\rm{(#1)}]}

\newcommand{\Res}{\mathop{\rm Res}\nolimits}

\newcommand{\scA}{{\mathcal A}^\circ}
\newcommand{\sB}{B^\circ}
\newcommand{\sG}{G^\circ}
\newcommand{\sM}{M^\circ}

\newcommand{\sQ}{Q^\circ}
\newcommand{\sS}{S^\circ}
\newcommand{\sU}{U^\circ}
\newcommand{\sX}{X^\circ}
\newcommand{\spi}{\bar\Pi(p)}

\newcommand{\pcA}{{\mathcal A}^\pi}
\newcommand{\pB}{B^\pi}

\newcommand{\pM}{M^\pi}

\newcommand{\pQ}{Q^\pi}
\newcommand{\pS}{S^\pi}
\newcommand{\pU}{U^\pi}
\newcommand{\pX}{X^\pi}

\renewcommand{\hat}{\widehat}

\newcommand{\mn}{\medskip\noindent}

\newtheorem{thm}{Theorem}[section]
\newtheorem{prop}[thm]{Proposition}
\newtheorem{cor}[thm]{Corollary}
\newtheorem{lem}[thm]{Lemma}
\newtheorem{conj}[thm]{Conjecture}
\newtheorem{exa}[thm]{Example}

\renewcommand{\qed}{\hfill $\square$}

\newcommand{\Definition}{\ \newline\refstepcounter{thm}\noindent{\bf Definition \thethm}\quad}

\newcommand{\Example}{\refstepcounter{thm}\noindent{\bf Example \thethm}\quad}
\newcommand{\Note}{\refstepcounter{thm}\noindent{\bf Note \thethm}\quad}

\newcounter{romanlistctr}
{\end{list}}%

 \makeatletter
 \@addtoreset{equation}{section}
 \makeatother

 \makeatletter
 \def\section{\@startsection {section}{1}{\z@}{-1.5ex plus -.5ex
 minus -.2ex}{1ex plus .2ex}{\large\bf}}

 \makeatletter
 \def\subsection{\@startsection {subsection}{1}{\z@}{-1.5ex plus -.5ex
 minus -.2ex}{1ex plus .2ex}{\bf}}

\begin{document}
\pagestyle{empty}
\title{A Quasi Curtis-Tits-Phan theorem for the symplectic group}
\author{Rieuwert J. Blok $^1$ and Corneliu Hoffman$^{1,2}$\\
$^1$Department of Mathematics and Statistics
\\ Bowling Green State
University \\ Bowling Green, OH 43403-1874\\[5pt]
$^2$ School of Mathematics \\
University of Birmingham \\
Edgbaston, B15 2TT, United Kingdom \\
}

\date{May 2008 \\[1in]
    \begin{flushleft}
    \end{flushleft}
       }

\maketitle

\newpage

\begin{flushleft} Proposed running head: \end{flushleft}
    \begin{center}
A quasi Phan-Curtis-Tits theorem for the symplectic group
    \end{center}

Send proofs to:
\begin{center}
    Rieuwert J. Blok\\ Department of Mathematics and Statistics\\
    Bowling Green State
University \\ Bowling Green, OH 43403\\[5pt]
Tel.: 419-372-7455\\
FAX: 419-372-6092\\
Email: blokr@member.ams.org
\end{center}
\newpage
\begin{abstract}
We obtain the symplectic group $\SP(V)$ as the universal completion of an amalgam of low rank subgroups akin to Levi components.
We let $\SP(V)$ act flag-transitively on the geometry of maximal rank subspaces of $V$.
We show that this geometry and its rank $\ge 3$ residues are simply connected with few exceptions.
The main exceptional residue is described in some detail.
The amalgamation result is then obtained by applying Tits' lemma.
This provides a new way of recognizing the symplectic groups from a small collection of small subgroups.
\end{abstract}

\newpage

\pagestyle{plain}

\section{Introduction}
In the revision of the classification  of finite simple groups one of the important steps requires one to
prove that if  a simple group $G$ (the minimal counterexample)  contains a certain amalgam of subgroups that
one normally finds in a known simple group $H$ then $G$ is isomorphic to $H$.
A geometric approach to recognition theorems was initiated in~\cite{BGHS03,BS04} and named Curtis-Phan-Tits
 theory.
The present paper uses a natural generalization to this theory to provide a new recognition theorem for symplectic groups.

Let us outline the Curtis-Phan-Tits theory setup. For details see~\cite{BGHS03}.
We consider a group $G$ which is either semi-simple of Lie type or a Kac-Moody group.
Let $\cT=(B_+, B_-)$ be the associated twin-building.
We first define  a {\em flip} to be an involutory automorphism $\sigma$ of $\cT$ that interchanges the two halves,
preserves distances and codistances and takes at least one chamber to an opposite.
Given a flip $\sigma$, construct $\cC_\sigma $ as the chamber system whose chambers are the pairs  of opposite
chambers $(c, c^\sigma)$ of $\cT$.
Let $G_\sigma$ be the fixed subgroup under the $\sigma$-induced automorphism of $G$.
Whenever the geometry $\Gamma_\sigma$ is simply connected one obtains $G_\sigma$ as the universal completion of the
amalgam of maximal parabolics for the action of $G_\sigma$ on $\Gamma_\sigma$.

We now exhibit a limitation of the Curtis-Tits-Phan setup using the building $\cT$ of type $A_n$ associated to $G=\PSL_{n+1}(\FF)$ for some field $\FF$.
In this setting a flip $\sigma$ is induced by a polarity.
The objects of the geometry $\Gamma_\sigma$ are the non-degenerate subspaces with respect to this polarity.
The requirement that $\sigma$-invariant pairs of opposite chambers exist enforces that at least one $1$-space $p$
does not intersect its polar hyperplane: that is, $p$ is non-absolute.
Polarities are classified by the Birkhoff-Von Neumann theorem and those with non-absolute points correspond to
 orthogonal or unitary forms.
Thus Curtis-Phan-Tits theory applied to the $A_n$ building only yields amalgams for orthogonal and unitary groups.

Our setup deviates from the general Curtis-Phan-Tits setup in the following way.
We still start with the building of type $A_n$, but we relinquish the requirement that $\Gamma_\sigma$ should consist of pairs of opposite chambers in $\cT$.
Instead we consider pairs of {\em almost opposite} chambers $(c,c^\sigma)$ corresponding to each other under a
 symplectic polarity $\sigma$. Note that $\sigma$ is {\em not a flip} since symplectic polarities have no non-absolute
 points.
We shall construct a geometry $\Gamma$ similar to $\Gamma_\sigma$ and obtain a presentation of the symplectic group.

There are several reasons why this setup and our result are interesting.
Recently it has become clear~\cite{DeMu,BlHo} that in order to study
(simple-) connectedness of Curtis-Phan-Tits geometries related to non-spherical twin-buildings one seems forced to study
 complexes of pairs of almost opposite chambers corresponding to each other under some involutory automorphism.
Our result involves a detailed study of just such a complex.

The geometry $\Gamma_\sigma$, which presented itself naturally to the authors as a Quasi Curtis-Tits-Phan geometry,
 is not new to the literature~\cite{Cu94,Ha88,Ha89,G04b}.
This indicates that this geometry is worth investigating for its own sake.

The fact that the geometry $\Gamma_\sigma$ has higher rank than the usual $C_n$ building geometry
 will cause the parabolic subgroups in our amalgam presentation to be smaller than those appearing in
 the amalgam of rank $2$ Levi-components.
This is a benefit by itself, but it will yield,  in particular, an amalgam presentation for $\Sp_4(q)$,
 a group that in the Curtis-Phan-Tits approach would be of rank $2$ itself and not admit an amalgam.

\paragraph{The main result}
We shall now briefly describe the main result.
Let $V$ be a vector space of dimension $2n\ge 4$ over a field $\FF$ endowed with a symplectic form $\sform$ of maximal
rank. Let $I=\{1,2,\ldots,n\}$.  Moreover consider $\cH=\{e_i,f_i\}_{i\in I}$ a hyperbolic basis of $V$.
The group of linear automorpisms of $V$ preserving the form $\sform$ will be called the
{\dfn symplectic group of $V$} and is denoted $G=\SP(V)$.
We shall define a geometry $\Gamma$ on the subspaces of $V$ whose radical has dimension at most $1$.
This geometry is transversal, residually connected, and simply connected.
Moreover, $G$ acts flag-transitively on $\Gamma$.
This then gives a presentation of $G$ as the universal completion of the amalgam $\cA$ of its maximal parabolic
subgroups with respect to its action on $\Gamma$.
Induction allows us to replace the amalgam $\cA$ of maximal parabolic subgroups by the amalgam $\cA_{\le 2}$ of
parabolic subgroups of rank at most $2$. A refinement of the amalgam $\cA_{\le 2}$ then leads to the following setup.

%

Consider the amalgam
$\cA^\pi=\{\pM_i,\pS_j,\pM_{i\, k},\pS_{j\, l},\pQ_{i\, j}\}_{j,l\in I; i,k\in I-\{n\}}$, whose groups are characterized as
follows:
\begin{itemize}
\item $\pS_i$ is the stabilizer in $G$ of all elements of $\cH-\{e_i,f_i\}$ and the subspace $\langle e_i,f_i\rangle$,
\item $\pM_i$ is the stabilizer in $G$ of all elements of $\cH-\{f_i,f_{i+1}\}$ and the subspace $\langle e_i,f_i,e_{i+1},f_{i+1}\rangle$,
\item $\pM_{i\, j}=\langle \pM_i,\pM_j\rangle$, $\pS_{i\, j}=\langle \pS_i,\pS_j\rangle$, $\pQ_{i\, j}=\langle \pM_i,\pS_j\rangle$.
\end{itemize}
We describe these groups in more detail in Section~\ref{section:slim}.
In particular, in Lemma~\ref{lem:structure sPij sQij sSij} we show that the groups in this amalgam are very small.
For instance,
$$\begin{array}{ll}
\pS_i & \cong \Sp_2(\FF),\\
\pM_i&\cong \FF^3.\\
\end{array}$$
Our main result is the following.
\bth\label{thm:main theorem1}
If $|\FF|>2$ and $V$ is a symplectic space of dimension $2n\ge 4$ over $\FF$ then the symplectic group $\Sp(V)$ is the universal completion of the amalgam $\cA^\pi$.
Moreover for any field $\FF$ the symplectic group is the universal completion of the amalgam $\cA$ from
Section \ref{subsection:parabolics}
\eth

\bco\label{mainco:sp4}
The group $\Sp_4(\FF)$  is the universal completion of the following amalgam.
$\{\pS_1,\pS_2,\pM_1,\pS_{1\, 2},\pQ_{1\, 1},\pQ_{1\, 2}\}$.
Here, $\pS_{1\, 2}\cong\Sp_2(\FF)\times\Sp_2(\FF)$ and
 $\pQ_{1\, 1}\cong\pQ_{1\, 2}\cong \FF^3\rtimes \Sp_2(\FF)$.
\eco


\setlength{\unitlength}{3pt}
\bpi(30,10)(-40,0)
\thicklines{
\multiput(3,3)(20,0){3}{\circle*{6}}
\multiput(5,3)(20,0){2}{\line(1,0){16}}
\put(3,-3){\makebox(0,0){$\pS_1$}}
\put(13,7){\makebox(0,0){$\pQ_{1\, 1}$}}
\put(23,-3){\makebox(0,0){$\pM_1$}}
\put(33,7){\makebox(0,0){$\pQ_{1\, 2}$}}
\put(43,-3){\makebox(0,0){$\pS_2$}}
}\epi

\bigskip
\paragraph{Organization of this paper}

In Section~\ref{section:preliminaries} we review some basic notions on geometries and some relevant facts about
 symplectic spaces.
In Section~\ref{section:geometry} we introduce a geometry $\Gamma$ on the almost non-degenerate subspaces of $V$
with respect to some symplectic form of maximal rank and describe its residues.
We prove that this geometry is transversal and residually connected.
In Section~\ref{section:simple connectedness} we show that the geometry and all residues of rank at least $3$ are
simply connected with one exception.
We show that in the exceptional case the residue has a simply connected $2$-cover.
In Section~\ref{section:group action} we describe the flag-transitive action of $\Sp(V)$ on $\Gamma$ and its rank $3$ residues.
We describe the parabolic subgroups in some detail and prove that $\Sp(V)$ is the universal completion of the amalgam
 $\cA_{\le 2}$ of parabolics of rank at most $2$ for its action on $\Gamma$.
In Section~\ref{section:slim} we define a slim version of the amalgam $\cA_{\le 2}$ by removing most of the Borel subgroup (for the action under consideration) from each of its groups.
Finally, in Section~\ref{section:concrete amalgam} we define the amalgam $\cA^\pi$ and prove
 Theorem~\ref{thm:main theorem1}.
\newpage
\section{Preliminaries}\label{section:preliminaries}
\subsection{Geometries}
For our viewpoint on geometries we'll use the following definitions from Buekenhout~\cite{Bu1995a}.
\Definition~\label{dfn:Buekenhout geometries}
A {\dfn pre-geometry} over a {\dfn type set} $I$ is a triple $\Gamma=(\cO,\typ,\inc)$, where $\cO$ is a collection of {\dfn objects} or {\dfn elements}, $I$ is a set of {\dfn types}, $\inc$ is a binary symmetric and reflexive relation, called the {\dfn incidence relation} and
 $\typ\colon \cO\to I$  is a {\dfn type function}  such that whenever $X\inc Y$, then either $X=Y$ or $\typ(X)\ne\typ(Y)$.

The {\dfn rank} of the pre-geometry $\Gamma$ is the size of $\typ(\cO)$.
A {\dfn flag} $F$ is a (possibly empty) collection of pairwise incident objects.
Its {\dfn type} (resp. {\dfn cotype}) is $\typ(F)$ (resp. $I-\typ(F)$).
The {\dfn rank} of $F$ is $\rank(F)=|\typ(F)|$.
The {\dfn type} of $F$ is $\typ(F)=\{\typ(X)\mid X\in F\}$. A {\dfn chamber} is a flag $C$ of type $I$.

A pre-geometry $\Gamma$ is a {\dfn geometry} if $\typ(\cO)=I$ and if $\Gamma$ is {\dfn transversal}, that is, if any flag is contained in a chamber.

The {\dfn incidence graph} of the pre-geometry $\Gamma=(\cO,\typ,\inc)$ over $I$ is the graph $(\cO,\inc)$.
This is a multipartite graph whose parts are indexed by $I$.
We call $\Gamma$ {\dfn connected} if its incidence graph is connected.

The {\dfn residue} of a flag $F$ is the pre-geometry
$\Res_\Gamma(F)=(\cO_F,\typ|_{\cO_F},\inc|_{\cO_F})$ over $I-\typ(F)$ induced on the collection $\cO_F$ of all objects in $\cO-F$ incident to $F$.
We call $\Gamma$ {\dfn residually connected} if for every flag of rank at least $2$
 the corresponding residue is connected.

For a subset $K\sbe I$ the {\dfn $K$-shadow} of a flag $F$ is the collection of all $K$-flags incident to $F$.

We will mostly be working with connected, residually connected geometries over a set $I$.
All our geometries will have a string diagram.
This means that these geometries have a {\dfn diagram} as in Buekenhout~\cite{Bu1995a} and this diagram will look
like a string, or a simple path.
For convenience of the reader we will give the following equivalent ad-hoc definition here.

\Definition
We say that a geometry {\dfn has a string diagram} if there is a total ordering on its type set $I$ such that for any three types $i,j,k\in I$ with $i<j<k$ we have the following.
If $X,Y,Z$ are objects of type $i$,$j$, and $k$ respectively such that $X$ and $Z$ are incident with $Y$, then $X$ is incident with $Z$.

Note that if a geometry has a string diagram, then so does every residue.

After choosing such a total ordering on $I$, we shall call the objects whose type is minimal in $I$ {\dfn points} and those objects whose type $i$ is minimal in
 $I-\{i\}$ are called {\dfn lines}.

\subsection{Automorphism groups and amalgams}
\Definition An {\dfn automorphism group} $G$ of a pre-geometry $\Gamma$ is a group of permutations of the collection of objects that preserves type and incidence.
We call $G$ {\dfn flag-transitive} if for any $J\sbe I$, $G$ is transitive on the collection of $J$-flags.

Let $G$ be a flag-transitive group of automorphisms of a geometry $\Gamma$ over an index set $I$.
Fix a chamber $C$.
The {\dfn standard parabolic subgroup of type $J\sbe I$} is the stabilizer in $G$ of the residue of type $J$ on $C$.

\Definition
In this paper we shall use the following definition of an amalgam of groups.
Let $(\cB,\prec)$ be a meet-semilattice with minimal element $\zh$ in which every maximal chain has length $s$.
An {\em amalgam} over $(\cB,\prec)$ is a collection of groups $\cA=\{A_\beta\mid \beta\in \cB\}$   together with a system of homomorphisms
 $\Phi=\{\varphi_{\beta,\gamma}\colon A_\beta\to A_\gamma\mid \beta\prec\gamma\}$
satisfying $\varphi_{\gamma,\delta}\after \varphi_{\beta, \gamma}=\varphi_{\beta,\delta}$
 whenever $\beta\prec\gamma\prec\delta$.
The number $s$ is called the {\em rank} of $\cA$.

The {\em universal completion} or {\em amalgamated sum} of ${\cA}$ is then a group $\hat{{G}}$ whose elements are words in the elements of the groups in $\cA$ subject to the relations between the elements of $A_\beta$ for any $\beta\in\cB$ and in which for each $\beta\prec\gamma$ each $a\in A_\beta$ is identified with $\varphi_{\beta,\gamma}(a)\in A_\gamma$. We then have a homomorphism $\hat{\cdot}\ \colon \cA\to \hat{{G}}$.

We note that for the appropriate choice of $(\cB,\prec)$ this definition of an amalgam and universal completion coincides with those given in \cite{Se1980,Ti1986b}.
\medskip

\Note
\begin{itemize}
\AI{i}
For each $\beta\in\cB$ we have a homomorphism\hspace{.3em} $\hat{\cdot}\ \colon A_\beta\to \hat{A_\beta}\le \hat{G}$, which is surjective, but not necessarily injective.
\AI{ii}
For $\beta,\gamma\in\cB$ with $\beta\prec \gamma$ we have $\hat{A_\beta}\le \hat{A_\gamma}$.
\AI{iii}
For $\beta,\gamma\in\cB$ we have $\hat{A_{\beta\meet\gamma}}\le\hat{A_\beta}\cap \hat{A_\gamma}$, but we do not a priori assume equality here.
\end{itemize}

\Example Let $G$ be a group acting flag-transitively on a geometry $\Gamma$ over
 an index set $I$.
Let $C$ be a chamber and, for every subset $J\sbe I$ with $|J|\le 2$ let $P_J$ be the standard parabolic subgroup of type $J$ in $G$.
Then, for $M\sbe K$ we have the natural inclusion homomorphisms
 $\varphi_{M,K}\colon P_M\to P_K$.
Hence $\cA=\{P_J\mid J\sbe I,\ |J|\le 2\}$ is an amalgam over
 $\cB=\{J\sbe I\mid |J|\le 2\}$ where $M\prec K\iff M\sbs K$.
For the universal completion $\hat{G}$ of $\cA$ we clearly have a surjective homomorpism $\tau\colon\hat{G}\to G$.

\mn
\subsection{Simple connectedness and amalgams}
In order to introduce the main tool of this paper, namely Lemma~\ref{lem:tits lemma} we need the notions of (closed) paths, (universal covers), simple connectedness, and the fundamental group.

In~\cite{Ti1986b,Fo1966,Ba1980,Qu1978} these notions are introduced in the context of (the face poset of) a simplicial complex in such a way that many classical results, such as can be found in \cite{Sp1981} continue to hold.
In the present paper we use definitions geared towards geometries. They are equivalent to those for the (face poset of) the simplicial complex, called the {\em flag complex} consisting of all flags of $\Gamma$ ordered by inclusion.
For a more extensive treatment of related issues see e.g.~\cite{Ti1986,Pas94}.

Let $\Gamma$ be a connected geometry over the finite set $I$.
A {\dfn path of length $k$} is a path $x_0,\ldots,x_k$ in the incidence graph.
We do not allow repetitions, that is, $x_i\ne x_{i+1}$ for all $0\le i<k$.
A {\dfn cycle based at an element $x$} is a path $x_0,\ldots,x_k$ in which $x_0=x=x_k$.
Two paths $\gamma$ and $\delta$ are {\dfn homotopy equivalent} if one can be obtained from the other by inserting
 or eliminating cycles of length $2$ or $3$. We denote this by $\gamma\simeq \delta$.
The homotopy classes of cycles based at an element $x$ form a group under concatenation.
This group is called the {\em fundamental group of $\Gamma$ based at $x$} and is denoted $\Pi_1(\Gamma,x)$.
If $\Gamma$ is (path) connected, then the isomorphism type of this group does not depend on $x$ and we call this group
 simply the {\em fundamental group of $\Gamma$} and denote it $\Pi_1(\Gamma)$.
We call $\Gamma$ {\em simply connected} if $\Pi_1(\Gamma)$ is trivial.

Given $k\in \NN_\ge 1$, a {\dfn $k$-covering} is an incidence and type preserving map $\pi\colon\bar{\Gamma}\to \Gamma$, where $\bar{\Gamma}$ and $\Gamma$ are geometries such that:
\begin{itemize}
\AI{CO1} For any non-empty $J$-flag $F$ in $\Gamma$ the fiber $\pi^{-1}(F)$ consists of exactly $k$ distinct and disjoint $J$-flags.
\AI{CO2} Given a non-empty flag $F$ in $\Gamma$ and some flag $\bar{F}$ in $\bar{\Gamma}$ such that $\pi(\bar{F})=F$, then the restriction
 $\pi\colon \Res(\bar{F})\to\Res(F)$ is an isomorphism.
\end{itemize}
We call $\bar{\Gamma}$ a {\em $k$-cover} of $\Gamma$.

The key consequence of this definition in terms of proving several classical results on covers is that they possess the {\em unique path lifting property}:
\begin{itemize}
\AI{L} Let $x=x_0,\ldots,x_n$ be a path in $\Gamma$ and suppose $\bar{x}\in\bar{\Gamma}$ satisfies $\pi(\bar{x})=x$. Then, there is a unique path $\bar{x}=\bar{x_0},\ldots,\bar{x_n}$ in $\bar{\Gamma}$ such that $\pi{\bar{x_i}}=x_i$ for all $i=0,1,\ldots,n$.
\end{itemize}

We call a cover of $\Gamma$ {\em universal} if it is universal in the category of coverings of $\Gamma$ and
 morphisms.
\bth\label{thm:simply connected}
Given a connected geometry $\Gamma$ and a covering $\pi\colon \bar{\Gamma}\to \Gamma$.
Then, the space $\bar{\Gamma}$ is universal among all covers of $\Gamma$, if and only if $\bar{\Gamma}$ is simply connected.
\eth
By a general categorical argument, a universal cover for $\Gamma$ is unique up to isomorphism.
Therefore $\Gamma$ is its own universal cover if and only if $\Pi_1(\Gamma)$ is trivial.

\mn

We can use the unique path lifting property to let $\Pi_1(\Gamma)$ act on $\bar{\Gamma}$ as follows.
Let $\delta=x_0,\ldots,x_n$ and $\bar{x}$ be as in (L) and denote by $[\delta]$ the homotopy class of $\delta$.
Then we define $[\delta]\cdot \bar{x_0}=\bar{x_n}$. This action has the following property (See \cite[Chap.6, Lemmas 1-4]{Sp1981} and \cite[Lemma 12.2]{Pas94}).

\ble\label{lem:action on fiber}
Let $\Gamma$ be a geometry with universal cover $(\bar{\Gamma},\pi)$ and let $x$ be an object.
Then the fundamental group $\Pi_1(\Gamma)$ acts regularly on $\pi^{-1}(x)$.
\ele
In Section~\ref{section:simple connectedness} we show that for some space $\Gamma$ we have $\Pi_1(\Gamma)\cong\ZZ_2$.
In view of the fact that $\Pi_1(\Gamma)$ is independent of the chosen base point, by Lemma~\ref{lem:action on fiber}
 this is equivalent to saying that the universal cover of $\Gamma$ is a $2$-cover.

\mn

The following result, which will be referred to as Tits' Lemma, is a consequence of~\cite[Corollaire 1]{Ti1986b}.
\ble\label{lem:tits lemma}
Given a group $G$ acting flag-transitively on a geometry $\Gamma$. Fix a maximal flag $C$.
Then $G$ is the universal completion of the amalgam consisting of the standard maximal parabolic subgroups of $G$ with respect to $C$ if and only if $\Gamma$ is simply connected.
\ele

\mn

We conclude this subsection with some methods to show that a geometry $\Gamma$ is simply connected.
This requires showing that any cycle based at a given element $x$ is homotopy equivalent to a cycle of length $0$.
We call such a cycle {\em trivial} or {\em null-homotopic}.

\ble\label{lem:common object}
If $\Gamma$ is a geometry, then any cycle all of whose elements are incident to a given element $A$  is null-homotopic.
\ele
\pf
In this case, the cycle together with the object $A$ forms a cone so that the cycle is null-homotopic.
\qed

\mn


\ble\label{lem:point-line cycles only}
Let $\Gamma$ be a connected, residually connected geometry over a set $I$ and let $i,j\in I$ be distinct.
Then, every cycle based at $x$ is homotopy equivalent to a cycle consisting of objects of type $i$ and $j$ only.
\ele
\pf
Let $x_0,x_1,\ldots,x_n=x_0$ be a cycle.
By transversality we may assume that $x_0$ is an $i$-object.
We proceed by induction on the number $N$ of objects of type $k\ne i,j$.
If $N=0$ we're done.
Suppose $N\ge 1$ and let $l$ be minimal so that $X=x_l$ has type different from $i$ and $j$.
By transversality we may assume that $x_{l+1}$ is of type $i$ or $j$.
By an easy induction argument one can show that since $\Gamma$ is residually connected, the incidence graph induced on the collection $X_{\{i,j\}}$ of all objects of type $i$ and $j$ incident to $X$ is connected.
Let $x_{l-1}=y_0,\ldots,y_m=x_{l+1}$ be a path in $X_{\{i,j\}}$.
Then $x_{l-1},x_l,x_{l+1}$ and $y_0,\ldots,y_m$ are homotopy-equivalent by Lemma~\ref{lem:common object}
 since all these objects are incident to $X$.
Thus we can replace $x_0,x_1,\ldots,x_l,\ldots,x_n=x_0$ by the homotopic path
 $x_0,\ldots,x_{l-1}=y_0,y_1,\ldots,y_m=x_{l+1},\ldots,x_n=x_0$ which contains only $N-1$ objects of type
 different from $i,j$.
By induction we are done.
\qed

\newpage

\section{A geometry for the symplectic group}\label{section:geometry}
Let $V$ be a vector space of dimension $n$ over a field $\FF$ endowed with a symplectic form $\sform$ of maximal rank.
We first need some notation to describe how $\sform$ restricts to the various subspaces of $V$.
Let $\perp$ denote the orthogonality relation between subspaces of $V$ induced by $\sform$.
Thus, for $U,W\le V$ we have
$$U\perp W \iff  \sform(u,w)=0\mbox{  for all }u\in U, w\in W.$$
We write
$$U^\perp=\{ v\in V\mid \sform(u,v)=0\  \forall u\in U\}\le V.$$

The {\dfn radical} of a subspace $U$ is the subspace $\Rad(U)=U\cap U^\perp$.
The {\dfn rank} of $U$ is $\rank(U)=\dim U-\dim\Rad(U)$.
Note that since $\sform$ is symplectic, we have $\sform(v,v)=0$ for all $v\in V$ and so $V$  has no anisotropic
 part with respect to $\sform$.

Let $2r=\rank(U)$ and $d=\dim(U)-2r$. A {\dfn hyperbolic basis} for $U$ is a basis $\{e_i,f_j\mid 1\le i\le r+d, 1\le j\le r\}$ such that
\begin{itemize}
\AI{i}
for all $1\le i,j\le r$,
$$\begin{array}{rl}
 \sform(e_i,e_j)=\sform(f_i,f_j)&=0,\\
 \sform(e_i,f_j)&=\delta_{i j},\mbox{ and }\\
 \end{array}$$
\AI{ii}
 $\{e_{r+i}\mid 1\le i\le d\}$ is a basis for $\Rad(U)$.
\end{itemize}

\ble\label{lem:hyperbolic basis}
Suppose that $W\le U-\Rad(U)$.
Then, any hyperbolic basis for $W$ extends to a hyperbolic basis for $U$.
\ele
\pf
This is in some sense Witt's theorem. See e.g.\ Taylor~\cite{Ta1992} for a proof.
\qed

\ble\label{lem:radical complement}
Suppose that $W\le U-\Rad(U)$ and $r=\rank(U)$.
If $\dim(W)=2r$, then $U=W\oplus \Rad(U)$ and $W$ is non-degenerate.
\ele

\pf
That $U=W\oplus \Rad(U)$ is simple linear algebra.
As a consequence, and since $\Rad(U)\perp U$, we have $\Rad(W)\le \Rad(U)$.
Thus $\Rad(W)\le \Rad(U)\cap W=\{0\}$.
\qed

\mn
\paragraph{The Quasi Curtis-Phan-Tits geometry}

The {\dfn quasi-Phan geometry} $\Gamma=\Gamma(V)$ is a geometry over $I=\{0,1,\ldots,n-1\}$ defined as follows.
For $i\in I$, the $i$-objects, or objects of type $i$, are the $i$-spaces $U\le V-\Rad(V)$ such that
 $\dim(\Rad(U))\le 1$.
More explicitly, since $\sform$ is symplectic, this means that
$$\dim(\Rad(U))=\left\{
\begin{array}{rl}
0 &\mbox{ if $i$ is even}\\
1 &\mbox{ if $i$ is odd.}\\
\end{array}\right.$$

We say that two objects $X$ and $Y$ are {\dfn incident} whenever $X\sbe Y-\Rad(Y)$ or vice versa.

\Note
The geometry $\Gamma$ defined above can also be defined simply using points and (hyperbolic) lines.
The objects of type $i\ge 3$ are then generated by $i$ points placed in suitable configuration
 (See for instance Lemma 3.2 in Blok~\cite{Bl2007}).
As such these geometries have been studied by Hall~\cite{Ha88,Ha89}, Cuypers~\cite{Cu94}, and
 Gramlich~\cite{G04b}.

\mn

It is not too difficult to see that $\Sp(V)$ is an automorphism group for $\Gamma$. In fact,
 Lemma~\ref{lem:flag-transitive} shows that it acts flag-transitively on $\Gamma$.

\bco\label{cor:hyperbolic bases and chambers}
Any hyperbolic basis for $V$ gives rise to a unique chamber of $\Gamma$ and, conversely, any chamber gives rise to a
 (not necessarily unique) hyperbolic basis for $V$.
\eco
\pf
Let $C$ be a chamber.
Then, for any two consecutive objects $W,U\in C\cup\{V\}$ we have $W\le U-\Rad(U)$.
Hence, using Lemma~\ref{lem:hyperbolic basis} repeatedly we find a hyperbolic basis $\cH$ for $V$ such that, for any $X\in C$, $\cH\cap X$ is a hyperbolic basis for $X$.

Conversely, let $\cH$ be a hyperbolic basis for $V$. Set $d=\dim(\Rad(V))$.
We have $\cH=\{e_i,f_j\mid 1\le i\le r+d, 1\le j\le r\}$.
For $1\le i\le r+d$, let $h_{2i-1}=e_i$ and, for $1\le j\le r$, let
$h_{2j}=f_j$.
Then, setting $C_l=\langle h_1,h_2,\ldots,h_l\rangle_V$, the collection $C=\{C_l\}_{l=1}^{n-1}$ is a chamber of $\Gamma$.
\qed

\medskip

In the remainder of this paper we shall use the following standard setup.
Let $2r=\rank(V)$  and let $d=\dim(\Rad(V))\le 1$, so that $n=2r+d$.
We fix a hyperbolic basis
$$\cH=\{e_i,f_j\mid 1\le i\le r+d, 1\le j\le r\}.$$
Alternatively we write
 $$\begin{array}{ll}
 \cH&=\{h_k\}_{k=1}^n\mbox{, where }\\
h_k&=\left\{\begin{array}{ll}
e_i&\mbox{ if }k={2i-1},\\
f_j & \mbox{ if }k=2j.\\
\end{array}\right.
\end{array}
$$
We call this the {\dfn standard hyperbolic basis}.
The {\dfn standard chamber} is the chamber $C=\{C_k\}_{k=1}^{n-1}$ associated to $\cH$ as in Corollary~\ref{cor:hyperbolic bases and chambers}.
That is, $C_k=\langle h_1,h_2,\ldots,h_k\rangle_V$, for all $1\le k\le n-1$.

\ble\label{lem:geometry transversal}
The pre-geometry $\Gamma$ is transversal and has a string diagram.
\ele
\pf
Let $F$ be a flag.
Then, for any two consecutive objects $W,U\in (F\cup\{V\})$ we have $W\le U-\Rad(U)$.
Hence, using Lemma~\ref{lem:hyperbolic basis} repeatedly we find a hyperbolic basis $\cH$ for $V$ such that, for any $X\in C$, $\cH\cap X$ is a hyperbolic basis for $X$.
According to Corollary~\ref{cor:hyperbolic bases and chambers}, $\cH$ defines a unique chamber $C$.
One verifies that $F\sbe C$.
The natural ordering on $I$ provides $\Gamma$ with a string diagram.
\qed

\ble\label{lem:geometry connected}
The pre-geometry $\Gamma$ is connected.
More precisely, the collinearity graph of the $\{1,2\}$-shadow geometry has diameter at most $2$ with equality if $n\ge 3$.
\ele
\pf
Let $X,Z$ be $1$-spaces in $V-\Rad(V)$.
If $X$ and $Z$ span a non-degenerate $2$-space, then we are done.
In particular, if $n=2$, then the diameter is $1$.

Other wise let $W$ be a point on $\langle X,Z\rangle$ different from $X$ and $Z$.
In case $\langle X,Z\rangle\spe \Rad(V)$, let $W=\Rad(V)$.
Since $V/\Rad(V)$ is non-degenerate, there is a point $Y$ in $W^\perp$ that is not
 in $\langle X,Z\rangle^\perp$.
Then clearly $X,Z\not\sbe Y^\perp$ and so $X,Y,Z$ is a path in $\Gamma$ from $X$ to $Z$.
Thus the collinearity graph of the $\{1,2\}$-shadow geometry has diameter at most $2$. Clearly equality holds.
\qed

\subsection{Residual geometries}
Let $C=\{C_i\}_{i\in I}$ be the standard chamber of $\Gamma$ associated to the hyperbolic basis $\cH$.
For every $J\sbe I$ the {\dfn standard residue of type $J$}, denoted $R_J$, is the residue of the $(I-J)$-flag
$\{C_i\}_{i\in I-J}$.
Let $\biguplus_{m=1}^M J_m$ be the partition of $J$ into maximal contiguous subsets.
(We call $K\sbe I$ contiguous if, whenever $i,k\in K$ and $i<j<k$, then $j\in K$.)

In this case, $R_J=R_{J_1}\times R_{J_2}\times\cdots\times R_{J_M}$ since $\Gamma$ has a string diagram.
It now suffices to describe $R_J$, where $J$ is contiguous.
Let $a=\min J$ and let $b=\max J$.
There are two cases according as $a$ is even or odd.

For odd $a$, the residue is the geometry $\Gamma(C_{b+1}/C_{a-1})\cong \Gamma((C_{a-1}^\perp\cap C_{b+1})/C_{a-1})$ of rank $b-a+2$. We set $C_0=\{0\}$ and $C_n=V$ for convenience.

For even $a$, we may assume that $V=(C_{a-2}^\perp \cap C_{b+1})/C_{a-2}$ and that $a=2$.
Thus we need to describe the residue of $C_1$.
We will show that $\Res_\Gamma(C_1)$ is isomorphic to a geometry $\Pi(p,H)$ defined as follows.

\Definition\label{dfn:PipH}
Note that $\dim(\Rad(V))\le 1$.
Let $p$ be a $1$-dimensional subspace of $V-\Rad(V)$ and let $H$ be some complement of $p$ in $V$ containing $\Rad(V)$.
Note that if $\Rad(H)$ is not trivial, then it is not contained in $p^\perp$.
Namely, $\Rad(H)\cap p^\perp\sbe \Rad(V)$, which is $0$ if $\dim(V)$ is even.

Then we define $\Pi(p,H)$ to be the geometry on the following collection of subspaces of $H$:
$$\{U\le H\mid \Rad(V)\not \subseteq U, \dim(\Rad(U))\le 2 \mbox{ and } \Rad(U) =\{0\}  \mbox{ or } \Rad(U)\not \subseteq  p^\perp\}.$$
Let $U$ and $W$ be in $\Pi(p,H)$ with $\dim U < \dim W$. We say that $U$ is incident to $W$ if either $\dim(W)$ is odd and $U\sbe W$ or $\dim(W)$ is even and $U\sbe W-\Rad(W\cap p^\perp)$.

More precisely, the objects are subspaces $U \le H$ not containing $\Rad(V)$ with the following properties
\begin{enumerate}
\item If $U$ is odd dimensional then $\Rad(U)$ has dimension $1$ and does not lie in $p^\perp$.
\item If $U$ is even dimensional then $U$ is either non-degenerate or $\Rad(U)$ has dimension $2$ and is not contained in $p^\perp$.
\end{enumerate}

\ble The map  $\varphi\colon\Res_\Gamma(p)\to\Pi(p,H)$ given by
$$X\mapsto X\cap H$$ is an isomorphism.
\ele
\pf
We first show that $X\in\Res_\Gamma(p)$ if and only if $X\cap H\in \Pi(p,H)$.
Note that since $p$ is isotropic, $p^\perp\cap H$ is a codimension $1$ subspace of $H$.
Note that if $X\in\Res_\Gamma(p)$, then since $X\in \Gamma$, $\Rad(V) \not \subset X$. Moreover  $X\cap H$ is a complement of $p$ in $X$ so $\Rad(X\cap H)$ cannot have  dimension more than two. Also $\Rad(X\cap H) \cap p^\perp \le \Rad(X)$ and so if $X$ is even dimensional  $\Rad(X\cap H) \cap p^\perp=\{0\}$ and if $X$ is odd dimensional and $X\cap H$ is degenerate then
 $\Rad(X\cap H) \not\sbe p^\perp$.

 Conversely if $U \in \Pi(p,H)$ it is easy to see that $X=\langle U,p\rangle$ is in $\Res_\Gamma(p)$. Indeed if $U$ is odd dimensional, since $p$ is not perpendicular to $\Rad(U)$, the space $X$ is non-degenerate. If $U$ is even dimensional and non-degenerate then the space $X$ is of maximal possible rank, $p$ is not in $\Rad(X)$ and $\Rad(V)\not\subseteq X$. Finally, if $U$ is even dimensional and $\Rad(U)$ has dimension $2$, then, since $p$ is not perpendicular to the whole of $\Rad(U)$, we have $\Rad(X)=\Rad(U)\cap p^\perp\ne \langle \Rad(V), p\rangle$ and so $X \in \Res_\Gamma(p)$.

 Suppose that $X$ and $Y$ are incident elements of $\Res_\Gamma(p)$ and $\dim(X)<\dim(Y)$. If $\dim(Y)$ is even then incidence is containment in both geometries. If $\dim(Y)$ is odd then we need to prove that $\Rad(Y)\not\sbe X$ iff $\Rad(Y^\varphi\cap p^\perp)\not\sbe X^\varphi$. We note that for any $Z\in \Res_\Gamma(p)$, we have
  $\Rad(Z)^\varphi=\Rad(Z^\varphi\cap p^\perp)$, and so the conclusion follows. \qed

 \ble\label{lem:odd point-line residual is complete graph}
 Two points $p_1,p_2$ of $\Pi(p,H)$ are collinear iff $\Rad(V) \not\sbe\langle p_1,p_2\rangle$. In particular if $\dim(V)$ is even then the collinearity graph of $\Pi(p,H)$ is a complete graph and if $\dim(V)$ is odd then the collinearity graph has diameter two.
 \ele
 \pf
If $p_1, p_2$ are two points in $\Pi(p,H)$, then $\langle p_1,p_2\rangle$ is either totally isotropic but not contained in $p^\perp$ or non-degenerate.
So this is a line of $\Pi(p,H)$ if and only if $\Rad(V) \not\sbe \langle p_1,p_2\rangle$.
Therefore the conclusion follows.
 \qed

\ble\label{lem:geometry residually connected}
The pre-geometry $\Gamma$ is residually connected.
\ele
Let $J\sbe I$ and let $R_J$ be the residue of the $(I-J)$-flag $\{C_i\}_{i\in I-J}$.
Let $\biguplus_{m=1}^M J_m$ be the partition of $J$ into contiguous subsets.
If $M\ge 2$, then the residue is connected since it is a direct product of geometries.
Otherwise, the residue is isomorphic to $\Gamma(V)$ for some vector space $V$ of dimension at least $3$, or to
 $\Pi(C_1,H)$ inside some $\Gamma(V)$ for some vector space $V$ of dimension at least $4$.
Thus the connectedness follows from Lemmas~\ref{lem:geometry connected} and~\ref{lem:odd point-line residual is complete graph}.
\qed

\bco\label{cor:geometry}
The pre-geometry $\Gamma$ is a geometry with a string diagram.
\eco
\pf
By Lemma~\ref{lem:geometry transversal}, $\Gamma$ is transversal and has a string diagram and by Lemma~\ref{lem:geometry residually connected}, it is residually connected.
\qed

\section{Simple connectedness}\label{section:simple connectedness}
In this section we prove that the geometry $\Gamma$ and all of its residues of rank at least $3$ are simply connected.
Note that $\Gamma$ and all its residues are geometries with a string diagram.
Therefore by Lemma~\ref{lem:point-line cycles only}  it suffices to show that all point-line cycles are null-homotopic.

\ble\label{lem:point-line cycles are null-homotopic}
If $|\FF|\ge 3$ or $\dim(V)$ is even, then any point-line cycle of $\Gamma$ is null-homotopic.
\ele
\pf
Let $\gamma$ be a point-line cycle based at a point $x$.
We identify $\gamma$ with the sequence $x_0,\ldots,x_k$ of points on $\gamma$ (so $\gamma$ is in fact a
$2k$-cycle).
We show by induction on $k$ that $\gamma$ is null-homotopic.

If $k\le 3$, then $U=\langle x_0,x_1,x_2\rangle_V$ is non-isotropic as it contains the hyperbolic line
$\langle x_0,x_1\rangle_V$.
If $\Rad(V)\not\sbe U$, then it is an object of the geometry and so by Lemma~\ref{lem:common object}, $\gamma$ is null-homotopic.
If $\Rad(V)\sbe U$, then since $|\FF|\ge 3$ there is a point $x$ such that
 $x$ is collinear to $x_0$, $x_1$, and $x_2$. Namely, consider the points
  $\Rad(V)x_i$ in the non-degenerate space $V/\Rad(V)$. These are $3$ distinct points on the non-degenerate $2$-space $U/\Rad(V)$. Take a fourth point $y$ on this line. Since $V/\Rad(V)$ is non-degenerate, there is a point $x$ orthogonal to $y$ but not to any other point on $U/\Rad(V)$.

Now let $k\ge 4$.
If two non-consecutive points $x_i$ and $x_j$ in $\gamma$ are collinear, then let
$\gamma_1=\delta_1 \circ \delta_2\circ\delta_1^{-1}$ and $\gamma_2=\delta_1\circ\delta_3$, where
 $$\begin{array}{rl}
 \delta_1&=x_0,\ldots,x_i\\
 \delta_2&=x_i,x_{i+1},\ldots,x_{j-1},x_j,x_i.\\
 \delta_3&=x_i,x_j,\ldots,x_k.\\
 \end{array}$$
Clearly $\gamma_1\circ\gamma_2$ is homotopic to $\gamma$.
Also, the cycles $\delta_2$ and $\gamma_2$ are both shorter than $\gamma$.
By induction, $\delta_2$ and $\gamma_2$ are null-homotopic and hence so are $\gamma_1$ and $\gamma$ itself.

Therefore we may assume that no two non-consecutive elements in $\gamma$ are collinear.
In particular $x_1\perp x_{k-1}$ and $x_{k-2}\perp x_0$.
The line $\langle x_{k-1},x_k\rangle_V$ has at least three points.
Take any $y\in\langle x_{k-1},x_k\rangle-\{x_{k-1},x_k\}$.
Then, $y\not\perp x_1,x_{k-2}$.
Thus, since $x_1\not\perp x_k=x_0$ and $x_1\perp x_{k-1}$, the point $y$ is collinear to $x_1$.
By the same reasoning $y$ is collinear to $x_{k-2}$.

Now let $\gamma_1=\delta_1\circ\delta_2$ and $\gamma_2=\delta_2^{-1}\circ\delta_3\circ\delta_2$ and $\gamma_3=\delta_2^{-1}\circ\delta_4\circ\delta_2$, where
 $$\begin{array}{rl}
 \delta_1&=x_0,x_1,y\\
 \delta_2&=y,x_0\\
 \delta_3&=y,x_1,x_2,\ldots,x_{k-2},y\\
 \delta_4&=y,x_{k-2},x_{k-1},y\\
 \end{array}$$
Then,
$$\gamma\simeq \gamma_1\circ\gamma_2\circ\gamma_3\simeq \gamma_2\simeq 0,$$
where the second equivalence holds since $\gamma_1$ and $\gamma_3$ are triangles, and the third equivalence holds since $\gamma_2\simeq 0$ by induction.
\qed

\bpr\label{prop:gamma simply connected}
The geometry $\Gamma$ is simply connected.
\epr
\pf
First of all, $\Gamma$ is connected by Lemma~\ref{lem:geometry connected}.
Thus it suffices to show that any cycle is null-homotopic.
By Lemma~\ref{lem:point-line cycles only} such a cycle is homotopic to a point-line cycle.
Finally, by Lemma~\ref{lem:point-line cycles are null-homotopic} point-line cycles are null-homotopic and the result follows.
\qed

\subsection{Residual geometries}
 \medskip
We shall now prove that, apart from one exception, all residues of rank at least $3$ are simply connected.
As we saw above any residue is either isomorphic to $\Gamma(V)$ for some $V$, or to $\Pi(p,H)$ for some point $p$ inside $\Gamma(V)$ for some $V$ or it is a direct product of such geometries and possibly rank $1$ residues.

We already proved that $\Gamma(V)$ is simply connected. We shall now prove that, apart from one exception, $\Pi(p,H)$ is a simply connected geometry. By Lemma~\ref{lem:odd point-line residual is complete graph}, if $n$ is even, then the collinearity graph of $\Pi(p,H)$ is a complete graph, and so we only have to show that all triangles are null-homotopic.
 If $n$ is odd, then the diameter of the collinearity graph is $2$ and so all $k$-cycles can be decomposed into cycles of length $3$, $4$, and $5$.
 \ble Suppose that $\dim(V)$ is odd. Then any cycle of length $4$ or $5$ can be decomposed into triangles.
 \ele
 \pf
We first note that two points $p_1$ and $p_2$ are at distance $2$ only if $\langle p_1,p_2\rangle\spe \Rad(V)$. Therefore, if $p_1,p_2,p_3,p_4$ is a $4$-cycle that cannot be decomposed into two triangles, then $\Rad(V)\sbe \langle p_1, p_3\rangle, \langle p_2,p_4\rangle$. Since $|\FF|\ge 2$ there exists a line $L$ on $\Rad(V)$ different from these two lines and any point of $L-p^\perp$ is collinear to all of $p_1,p_2,p_3,p_4$. Thus we decompose the $4$-cycle into triangles.
Now suppose $p_1,p_2,p_3,p_4,p_5$ is a  $5$-cycle. Then $\Rad(V)$ lies on at most $1$ of $p_1p_3$ and $p_1p_4$ and so one of these lines is in fact a line of the geometry. Thus, we can decompose the $5$-cycle into shorter cycles.
 \qed

 \ble\label{lem:hyperplanes of Pi} Consider a hyperplane $W$ of $H$.
 If $\dim(V)$ is even, then $W\in\Pi(p,H)$.
 If $\dim(V)$ is odd, then $W\in\Pi(p,H)$ if and only if $\Rad(V)\not\sbe W$.
 \ele
 \pf
 Let $\dim(V)$ be even, and let $S=\Rad(H)$. Then $\dim(S)=1$ and $S\not\sbe p^\perp$.
 If $S\not\sbe W$, then $W$ is non-degenerate, and so it belongs to $\Pi(p,H)$.
 If $S\sbe W$, then $S\sbe \Rad(W)$, which has dimension $2$ and so $\Rad(W)\not\sbe p^\perp$.
 So again $W\in \Pi(p,H)$.

 Now let $\dim(V)$ be odd, and let $R=\Rad(V)$. Then $R\sbe \Rad(H)$, which has dimension $2$ and is not included in $p^\perp$. Now since $W$ is a hyperplane of $H$ it either contains $\Rad(H)$ or it intersects it in a $1$-dimensional space. In the former case, $R\sbe \Rad(W)$ and $W\not\in\Pi(p,H)$.
 In the latter case, either $R\sbe \Rad(W)$ and $W\not\in\Pi(p,H)$, or $\Rad(W)\not\sbe p^\perp$ and
  $W\in\Pi(p,H)$.
 \qed

 \ble\label{lem:Pi simply connected n odd}
If $\dim(V)$ is odd, then any triangle of $\Pi(p,H)$  is null-homotopic.
 \ele
 \pf
 Take a triangle on the points $p_1,p_2,p_3$. Note that this means that $R=\Rad(V)$ does not lie on any of the lines $p_ip_j$. Let $U=\langle p_1,p_2,p_3\rangle$.
 We have two cases:
 1) $R\not\sbe U$. Then pick a hyperplane $W$ of $H$ containing $U$ but not $R$.
 Then by Lemma~\ref{lem:hyperplanes of Pi}, $W\in \Pi(p,H)$. Moreover, since $\dim(W)$ is odd it is incident to each of the lines $p_ip_j$.

 2) $R\sbe U$. Let $S$ be a point of $\Rad(H)-R$. Then $S\not\sbe p^\perp$.
 Also note that for any $i$ and $j$, $p_i,p_j,S$ is a triangle of type (1), so the triangle  $p_1,p_2,p_3$ is null-homotopic.
 \qed

 \ble\label{lem:>6 or >2 is null}
 Suppose $\dim(V)\ge 7$ or $|\FF|\ge 3$.
 Then any triangle of $\Pi(p,H)$  is null-homotopic.
 \ele
\pf
If $\dim(V)$ is odd, then we are done by Lemma~\ref{lem:Pi simply connected n odd}.
So now let $\dim(V)=n$ be even.
Let $Q=\Rad(H)$.
Take a triangle on the points $p_1,p_2,p_3$ and let $U=\langle p_1,p_2,p_3\rangle$.
If $U$ is a plane, then we're done by Lemma~\ref{lem:common object}.
From now on assume this is not the case. So $U$ is either totally isotropic, or it has rank $2$ and  $T=\Rad(U)$ is contained in $p^\perp\cap H$.

Let $L=U\cap p^\perp$. For $1\le i<j\le 3$, let $L_{i j}=p_ip_j$ and $q_{i j}=L_{i j}\cap L$.
Our aim is to find an object $W$ of $\Pi(p,H)$ that is incident to all lines $L_{i j}$.
Then, by Lemma~\ref{lem:common object} we are done.

Suppose we can find a point $r$ in $p^\perp\cap H$ such that $r^\perp\spe L$.
For such a point let $W=\langle r^\perp\cap p^\perp\cap H, U\rangle $.
First note that $\dim(r^\perp\cap p^\perp\cap H)=n-3$.
This is because $p^\perp\cap H$ is a non-degenerate symplectic space of even dimension $n-2$. Clearly $r^\perp\cap p^\perp\cap H$ has radical $r$.

Note that $L\sbe  r^\perp\cap p^\perp\cap H$.
Also, since $U\not\sbe p^\perp\cap H$ and $U\cap p^\perp = L$ we have
 $U\cap (W\cap p^\perp)=L$ and $W\cap p^\perp=r^\perp\cap p^\perp\cap H$.
Hence $W=\langle (W\cap p^\perp),U\rangle$ and so $\dim(W)=\dim(W\cap p^\perp)+\dim(U)-\dim(L)=n-2$.
Thus $W$ is an object of $\Pi(p,H)$ by Lemma~\ref{lem:hyperplanes of Pi} as required.

We now show that we can find such a point $r$ that is different from $q_{1 2}$ , $q_{1 3}$, and $q_{2 3}$. Then, since $r=\Rad(W\cap p^\perp)$, all lines $L_{i j}$ are incident to $W$ as required.

There are two cases.
(1) $n=6$ and $|\FF|\ge 3$.
Since $|\FF|\ge 3$ we can find a point $r\in L-\{q_{1 2},q_{1 3},q_{2 3}\}$.
Note that $L$ is totally isotropic and so $L\sbe r^\perp$, as required.

(2) $n\ge 7$. Note that in this case in fact $n\ge 8$ since we assume that $n$ is even.
We find a point $r$ in $L^\perp\cap p^\perp\cap H-L$.
This is possible since $p^\perp\cap H$ is non-degenerate and $L$ is a $2$-space so that $\dim(L^\perp\cap p^\perp\cap H)=(n-2)-2\ge 4>\dim(L)=2$. Again we have found $r$ as required.
\qed

\bpr\label{prop:residues simply connected}
If $|\FF|\ge 3$ then the residues of rank at least $3$ are simply connected. \epr
\pf
Let $J\sbe I$ and let $R_{I-J}$ be the residue of the $J$-flag $\{C_j\}_{j\in J}$.
Set $r=|I-J|$.
Let $\biguplus_{m=1}^M I_m$ be the partition of $I-J$ into contiguous subsets.
If $M\ge 2$, then the residue is a direct product of two geometries, at least one of which is a residue of rank at least $2$. Such rank $2$ residues are connected by Lemma~\ref{lem:geometry residually connected}.
Therefore $R$ is simply connected in this case.
Otherwise, the residue is isomorphic to $\Gamma(V)$ for some vector space $V$ of dimension at least $3$, or to
 $\Pi(p,H)$ for some point $p$ inside some $\Gamma(V)$ for some vector space $V$ of dimension at least $4$.
Therefore the simple connectedness follows from Lemmas~\ref{lem:>6 or >2 is null},~\ref{lem:Pi simply connected n odd} and Proposition~\ref{prop:gamma simply connected}.
\qed

\medskip
\subsection{The exceptional residue}
We are now left with the intriguing case $n=6, |\FF|=2$. Let us first describe the geometry $\Pi(p,H)$.
The points are the $16$ points of $H-p^\perp$.
The lines are those lines of $H$ not contained in $p^\perp$.
Thus there are two types of lines, totally isotropic and hyperbolic ones. Each line has exactly two points and any two points are on exactly one line.

The planes of $\Pi(p,H)$ are those non-isotropic planes of $H$ whose radical is not contained in $p^\perp$.
Such a plane has rank $2$ and its radical is a point of $\Pi(p,H)$, i.e.\ it is in $H-p^\perp$.
Each plane has exactly $4$ points. A plane may or may not contain $\Rad(H)$. Any point or line contained in a plane is incident to that plane.

The $4$-spaces of $\Pi(p,H)$ are those $4$-spaces of $H$ that are either non-degenerate or have a radical of dimension $2$ that is not contained in $p^\perp$.
A $4$-space $W$ is incident to all $8$ points it contains and is incident to any line or plane it contains that doesn't pass through the point $\Rad(W\cap p^\perp)$.
Thus in fact $W$ is incident to all planes $Y$ of $\Pi(p,H)$ contained in $W$. Namely, if $\Rad(W\cap p^\perp)\le Y$, then $p^\perp\cap Y$ is a totally isotropic line, implying that $\Rad(Y)\le p^\perp$, a contradiction.

The geometry $\Pi(p,H)$ is not simply connected.
To see this, we construct a simply connected rank $4$ geometry $\spi$ which is a degree two cover of $\Pi(p,H)$.
Let us now describe the geometry $\spi$.
A point $q\in \Gamma$ is a point of $\spi$ if the line $pq$ is non-degenerate. In particular the points of $\Pi(p,H)$ are among the points of $\spi$. We construct  a map, which on the point set of $\spi$ is given by
$$\begin{array}{rl}
\psi\colon \spi &\to \Pi(p,H)\\
q &\mapsto\langle p,q\rangle\cap H.
\\
\end{array}$$
This map is two-to-one and for any $q\in \Pi(p,H)$ we set $q^-=q$ and denote by $q^+$ the point of $\spi$ such that  $\psi^{-1}(q)=\{q^-,q^+\}$. We extend this notation to arbitrary point sets $S$ of $\Pi(p,H)$ by setting $S^\epsilon=\{q^\epsilon\mid q\in S\}$ for $\epsilon=\pm$.

We now note that every object of $\Pi(p,H)$ can be identified with its point-shadow. It follows from the above description of $\Pi(p,H)$ however that inclusion of point-shadows does not always imply incidence.

In order to describe $\spi$, we shall identify objects in $\Pi(p,H)$ and $\spi$ with their point-shadow. We denote point-shadows by roman capitals. If we need to make the distinction between objects and their point-shadows explicit, we'll use calligraphic capitals for the objects and the related roman capitals for their point-shadows.
Thus $\cX$ may denote an object of $\Pi(p,H)$ (or $\spi$) whose point-shadow is $X$.

For any object $\cX$ of $\Pi(p,H)$ we will define exactly two objects $\cX_-$ and $\cX_+$ in $\spi$ such that $\psi(X_-)=\psi(X_+)=X$ as point-sets.
We then define objects $\cX$ and $\cY$ of $\spi$ to be incident whenever
 $X\sbe Y$ or $Y\sbe X$ and $\psi(\cX)$ and $\psi(\cY)$ are incident in $\Pi(p,H)$.

We shall obtain $\cX_-$ and $\cX_+$ by defining a partition $X_0\uplus X_1$ of the point-set of $\cX$ and setting
 $X_+=X_0^+\uplus X_1^-$ and $X_-=X_0^-\uplus X_1^+$.

First let $\cX$ be a line. If $\cX$ is non-degenerate, then $X_0=X$ and $X_1=\emptyset$ so that $X_+=X^+$ and $X_-=X^-$. If $X$ is isotropic, then $X=\{x_0\}$ and $X_1=\{x_1\}$ so that $X_-=\{x_0^-,x_1^+\}$ and $X_+=\{x_0^+,x_1^-\}$.

Next, let $\cX$ be a plane. Then $X_0=\{\Rad(X)\}$ and $X_1=X-X_0$.
Note that for every line $\cY$ incident to $\cX$ the partition $Y_0\uplus Y_1$ agrees with the partition $X_0\uplus X_1$.
Hence if $\cY$ is a line and $\cX$ is a plane in $\spi$ such that $\psi(\cY)$ is incident to $\psi(\cX)$, then either $Y\sbe X$ or $Y\cap X=\emptyset$.

Finally, let $\cX$ be a $4$-space. Let $r=\Rad(X\cap p^\perp)$. Note that the projective lines $L_i$ ($i=1,2,3,4$) of $\cX$ on $r$ meeting $H-p^\perp$ are not incident to $\cX$ in $\Pi(p,H)$.

First let $\cX$ be non-degenerate. Then $L_i$ is non-degenerate for all $i$.
For $i=1,2,3,4$, let $L_i=\{r,p_i,q_i\}$ such that $q_i=p_1^\perp\cap L_i$ for $i=2,3,4$. Note that this implies the fact that if $i\ne j$, then  $p_i\perp q_j$ but $q_i\not\perp q_j$ and $p_i\not\perp p_j$.
Then $X_0=\{p_1,\ldots,p_4\}$ and $X_1=\{q_1,\ldots,q_4\}$.
We now claim that for every line $\cY$ incident to $\cX$ the partition $Y_0\uplus Y_1$ agrees with $X_0\uplus X_1$.
Thus we must show that if $\cY$ is non-degenerate, then $Y\sbe X_0$ or $Y\sbe X_1$ and if $\cY$ is isotropic, then
 $Y$ intersects $X_0$ and $X_1$ non-trivially.
First, since $\cX$ is non-degenerate, the lines $q_iq_j$ are all non-degenerate for $2\le i<j\le 4$.
Moreover, $q_1q_i$ is non-degenerate as well, for $i=2,3,4$ since otherwise $L_i$ must be totally isotropic, a contradiction.  As a consequence, $q_1p_i$ is totally isotropic. Hence by interchanging the $p_i$'s for $q_i$'s we see that $p_ip_j$ is non-degenerate for all $1\le i<j\le 4$.
Considering the non-degenerate lines $L_i$ and $L_j$ we see that $p_iq_j$ must be isotropic for all $1\le i\ne j\le 4$.
We have exhausted all lines incident to $\cX$ and it is clear that the claim holds.

If $\cX$ is degenerate, then it has a radical $R=\Rad(X)$ of dimension $2$ passing through the radical $r=\Rad(X\cap p^\perp)$. In this case the lines $L_i$ are all totally isotropic. Let $R=L_1$. Then $X_0=\{p_1,q_1\}$ and $X_1=\{p_i,q_i\mid i=2,3,4\}$.
We now claim that for every line $\cY$ incident to $\cX$ the partition $Y_0\uplus Y_1$ agrees with $X_0\uplus X_1$.
Clearly every line meeting $X_0$ and $X_1$ is totally isotropic. Next consider a line $\cY$ meeting $L_i$ and $L_j$
 in $X_1$. If $\cY$ were totally isotropic, then $X=\langle L_1,L_i,L_j\rangle$ is totally isotropic, a contradiction.
We have exhausted all lines incident to $X$ and it is clear that the claim holds.

Next we claim that for every plane $\cY$ incident to $\cX$ the partition $Y_0\uplus Y_1$ agrees with $X_0\uplus X_1$.
First note that $\cY$ doesn't contain $\Rad(X\cap p^\perp)$, because that would make $Y\cap p^\perp$ totally isotropic and $\Rad(Y)\le p^\perp$ contrary to the description of planes of $\Pi(p,H)$.
As a consequence, if $\cL$ is a line incident to $\cY$, then it doesn't meet $\Rad(X\cap p^\perp)$ and so is incident to $\cX$ as well.
Hence since for every line $\cL$ incident to $\cY$, the partition $L_0\uplus L_1$ agrees with the partition $X_0\uplus X_1$ as well as with the partition $Y_0\uplus Y_1$, also the latter two partitions agree.
We can conclude that if $\cY$ is a point, line or plane and $\cX$ is a $4$-space of $\spi$ such that
 $\psi(\cY)$ and $\psi(\cX)$ are incident, then either $Y\sbe X$ or $Y\cap X=\emptyset$.

Let $\Psi$ be the extension of $\psi$ to the entire collection of objects in $\spi$.
\ble\label{lem:2-cover}
\begin{itemize}
\AI{a} The map $\Psi\colon\spi\to\Pi(p,H)$ is $2$-to-$1$ on the objects.
For any object $\cX$ of $\Pi(p,H)$, if $\Psi^{-1}(\cX)=\{\cX_-,\cX_+\}$, then $\psi^{-1}(X)=X_-\uplus X_+$.
\AI{b} Two objects $\cal{X}$ and $\cal{Y}$ of $\spi$ are incident if and only if $\Psi(\cal{X})$ and $\Psi(\cal{Y})$ are incident and $X\cap Y\ne\emptyset$.
\AI{c}
For any non-empty flag $F_\bullet$ in $\spi$, the map $\Psi\colon\Res(F_\bullet)\to\Res(\Psi(F_\bullet))$ is an isomorphism of geometries.
\AI{d} The map $\Psi\colon\spi\to\Pi(p,H)$ is a $2$-cover.
\AI{e} The pre-geometry $\spi$ is transversal.
\end{itemize}
\ele
\pf
(a) This clear from the construction of $\spi$.

(b) By definition of incidence in $\spi$, $\cX$ and $\cY$ can only be incident if $\Psi(\cX)$ and $\Psi(\cY)$ are incident. The preceding discussion has shown that in this case either $X\cap Y=\emptyset$ or $X\sbe Y$ or $Y\sbe X$. Now $\cX$ and $\cY$ are defined to be incident precisely in the latter case.

(c)
Let $\cX_\bullet$ and $\cY_\bullet$ denote objects of $\spi$.
Also, let $\cX=\Psi(\cX_\bullet)$, $\cY=\Psi(\cY_\bullet)$ and $F=\Psi(F_\bullet)$.
By definition of the objects in $\spi$, $\psi\colon X_\bullet\to X$ is a bijection of points.
Therefore if $\cY$ is incident with $\cX$, then $\cX_\bullet$ is incident with exactly one object in $\Psi^{-1}(\cY)$.
Let $\cY$ be incident with $F$. Then by the same token $F_\bullet$ is incident with at most one object in $\Psi^{-1}(\cY)$.
We now show that there is at least one such object.
Suppose that $\cY_\bullet$ is incident to at least one object $\cZ_\bullet$ of $F_\bullet$.
Without loss of generality assume $Y_\bullet\sbe Z_\bullet$.
If $X_\bullet$ is an element of $F_\bullet$ and $Z_\bullet\sbe X_\bullet$ then $\cY_\bullet$ is incident to $\cX_\bullet$ as well.
Now assume $X_\bullet \sbe Z_\bullet$.
Since $\cX$ and $\cY$ are incident, $\cY_\bullet$ must be either incident to $\cX_\bullet$ or to the object in $\Psi^{-1}(\cX)$ different from $\cX_\bullet$.
In the latter case it follows that $\cZ_\bullet$ is incident to both objects in $\Psi^{-1}(\cX)$, a contradiction.
A similar argument holds when $Z_\bullet\sbe Y_\bullet$.
We conclude that $\Psi\colon \Res(F_\bullet)\to \Res(F)$ is a bijection.
Clearly incidence is preserved by $\Psi$, but we must show the same holds for $\Psi^{-1}\colon \Res(F)\to\Res(F_\bullet)$. Let $\cX,\cY\in\Res(F)$ be incident and let $\cX_\bullet, \cY_\bullet\in \Res(F_\bullet)$.
Then there is a point $q$ incident to $\cX$, $\cY$ and $\cF$.
Suppose $q^\epsilon\in F_\bullet$. Then $q^\epsilon\in X_\bullet\cap Y_\bullet$ and by (b) we find that $\cX_\bullet$ and $\cY_\bullet$ are incident.

(d) This follows from (a) and (c).

(e) This is immediate from (c).
\qed

\ble\label{lem:spi connected}
The pre-geometry $\spi$ is  connected. Any two points are at distance at most $2$, except the points
 $\Rad(H)^\pm$, which are at distance $3$ from one another.
\ele
\pf
Let $\epsilon\in\{+,-\}$.
Let $Q=\Rad(H)$.
Then for any point $q\ne Q$, since the line $qQ$ is totally isotropic, $Q^\epsilon$ is collinear to $q^{-\epsilon}$ but not to $q^\epsilon$.
In particular any two points with the same sign are at distance at most $2$. It is also clear that $Q^+$ and $Q^-$ have no common neighbors and are at distance at least $3$.

Now consider two points $q_1,q_2\ne Q$. If the line $q_1q_2$ is totally isotropic, then $q_1^\epsilon$ is collinear to $q_2^{-\epsilon}$ in $\spi$.
If the line $q_1q_2$ is non-degenerate, we claim that there exists $q_3\in\Pi(p,H)$ with $q_1\perp q_3\not\perp q_2$.
Namely, we must show that $q_1^\perp-(q_2^\perp\cup p^\perp)\ne\emptyset$. However, this is clear since both
 $p^\perp$ and $q_2^\perp$ define proper hyperplanes of the $4$-space $q_1^\perp$. Since no linear subspace of $V$ is the union of two of its hyperplanes our claim follows.
In $\spi$ we find both $q_1^\epsilon$ and $q_2^{-\epsilon}$ collinear to $q_3^{-\epsilon}$ so again $q_1$ and $q_2$ are at distance at most $2$.

Finally consider $Q^+$ and $Q^-$. Let $q_1$ be a point of $\Pi(p,H)$ and let $q_2$ be a point of $\Pi(p,H)$
 in $q_1^\perp-\{Q\}$. Then $Q^-,q_1^+,q_2^-,Q^+$ is a path of length $3$.
\qed

\ble\label{lem:spi residually connected}
The  pre-geometry $\spi$  is  residually  connected.
\ele
\pf
By Lemma~\ref{lem:spi connected}, $\spi$ is connected so it suffices to show that every residue of rank at least $2$ is connected. This follows immediately from Lemmas~\ref{lem:2-cover}~and~\ref{lem:geometry residually connected}.
\qed

\ble\label{lem:spi geometry}
The pre-geometry $\spi$ is a geometry with a string diagram.
\ele
\pf
By Lemma~\ref{lem:2-cover} it is transversal and by Lemma~\ref{lem:spi residually connected} it is residually connected. Therefore it is a geometry. That it has a string diagram is clear since it is a cover of $\Pi(p,H)$, which does have a string diagram.
\qed

\ble\label{lem:spi simply connected}
The geometry $\spi$ is simply connected.
\ele
\pf
By Lemma~\ref{lem:spi geometry} and ~\ref{lem:point-line cycles only} it suffices to  show that any point-line cycle is null-homotopic.

Let $Q=\Rad(H)$. For any point $q\in \Pi(p,H)$, let $q^*$ denote one of $q^+$, $q^-$.
We claim that any $k$-cycle with $k\ge 5$ can be decomposed into triangles, quadrangles, and pentagons.
Namely, let $\gamma=q_1^*,q_2^*,\ldots,q_k^*,q_1^*$ be a $k$-cycle in $\spi$.
If $q_1$ and $q_4$ are not both $Q$, then they are at distance at most $2$ by Lemma~\ref{lem:spi connected} and so
 we can decompose $\gamma$ as $(q_1^*,q_2^*,q_3^*,q_4^*)\after \delta \after\delta^{-1}\after
  (q_4^*,\ldots,q_k^*,q_1^*)$, where $\delta$ is a path from $q_1^*$ to $q_4^*$ of length at most $2$.
Thus  we can decompose the $k$-cycle into a $(k-1)$-cycle and a quadrangle or pentagon.
If $q_1$ and $q_2$ are both $Q$, then replacing $q_1$ and $q_4$ by $q_2$ and $q_5$, we can again decompose the $k$-cycle into a $(k-1)$-cycle and a quadrangle or pentagon.

We shall now analyze the triangles, quadrangles, and pentagons case by case.

\paragraph{Triangles}
The points of a triangle $q_1^*,q_2^*,q_3^*$ either all have the same sign or one has a sign different from the others. Note that a point that is collinear to another point with the same sign can not  cover $Q$. In both cases $X=\langle q_1,q_2,q_3\rangle$ has dimension $3$.
Also, $X$ is non-isotropic because at least one of the lines $\langle q^*_i,q^*_j\rangle$ is non-isotropic.
We'll show that $r=\Rad(X)$ does not lie in $p^\perp$.
In the case where $q_1^*$, $q_2^*$, and $q_3^*$ all have the same sign, $q_1,q_2,q_3$ form a triangle in $\Pi(p,H)$ whose lines are non-degenerate. In particular $r$ does not lie on any of these lines. The three remaining points on these lines are $X\cap p^\perp$.
In the latter case this is because $r$ is covered by the point of the triangle with the deviating sign.
Thus $X$ is a plane of $\Pi(p,H)$ and $q_1^*,q_2^*,q_3^*$ belong to a plane of $\spi$.

\paragraph{Quadrangles}
Next we consider a quadrangle $q_1^*,q_2^*,q_3^*,q_4^*$.
There are four cases.
Let $\epsilon\in\{+,-\}$.

(1) All points have the same sign, say $\epsilon$. As we saw with the triangles, none of these is $Q^\epsilon$. Hence $Q^{-\epsilon}$ is connected to all $4$ points and so the quadrangle decomposes into triangles.

(2) All points but one, have the same sign, say $+$. Note that again the points on the same level do not cover $Q$, but are all collinear to $Q^-$. This decomposes the cycle into triangles and a quadrangle with two points of each sign.

(3) There are two points of each sign and these points are consecutive.
Without loss of generality let $q_1^-,q_2^-,q_3^+,q_4^+$ be the quadrangle.
We first note that the points $q_1,q_2,q_3,q_4$ are all distinct.
This is because no point of $\spi$ is collinear to both covers of the same point in
 $\Pi(p,H)$.
Also note that $q_1q_4$ and $q_2q_3$ are totally isotropic, but $q_1q_2$ and $q_3q_4$ are not. We may also assume that $q_1q_3$ and $q_2q_4$ are non-degenerate, for otherwise we can decompose the quadrangle into triangles.
Consider the space $Y=q_1^\perp\cap q_2^\perp\cap H$. It is a $3$-dimensional space  of rank $2$ whose radical is $Q$.
 Now both $q_3^\perp\cap Y$ and $q_4^\perp\cap Y$ are lines of $Y$ through $Q$.
 Note that $p^\perp\cap Y$ is a line of $Y$ not through $Q$.
Hence there is a point $q\in Y-q_3^\perp-q_4^\perp-p^\perp$.
We find that $q^+$ is collinear to all points on the quadrangle $q_1^-,q_2^-,q_3^+,q_4^+$, which therefore decomposes into triangles.

(4) There are two points of each sign and these points are not consecutive.
Without loss of generality let $q_1^-,q_2^+,q_3^-,q_4^+$ be the quadrangle.
Let $X=\langle q_1,q_2,q_3,q_4\rangle$.
We first note that we can assume that the points $q_1,q_2,q_3,q_4$ are all distinct. No two consecutive ones can be the same so if for example $q_1=q_3$ then the quadrangle is just a return.
 Hence $\dim(X)=3,4$.
Note that if either $q_1q_3$ or $q_2q_4$ is non-degenerate, then we can decompose the quadrangle into triangles. Therefore all lines $q_iq_j$ are totally isotropic and it follows that $X$ is a totally isotropic $3$-space.
This means that $Q=q_i$ for some $i$, which we may assume to be $4$.
Consider a totally isotropic $3$-space $Y$ on $q_2q_4$ different from $X$.
Then $q_2q_4$ and $p^\perp \cap Y$ are intersecting lines of $Y$ and so there is a point
 $q\in Y-q_2q_4-p^\perp$.
Note that $q\not \perp q_1,q_3$ for otherwise $\langle q, X\rangle$ is a totally isotropic $4$-space.
We find that $q^-$ is collinear to all points of the quadrangle  $q_1^-,q_2^+,q_3^-,q_4^+$, which therefore decomposes into triangles.

\paragraph{Pentagons}
We first note that we may assume that a pentagon has no more than $2$ consecutive points of the same sign.
If it contains $4$ or more of sign $\epsilon$, then these points are all collinear to $Q^{-\epsilon}$ which then yields a decomposition of the pentagon into triangles and quadrangles. If it contains  exactly $3$ consecutive points at the same level   the same argument decomposes it into $2$ triangles and a pentagon that contains no more than $2$ consecutive points at the same level.

Therefore we may assume without loss of generality that the pentagon is
 $q_1^-$, $q_2^-$, $q_3^+$, $q_4^-$, $q_5^+$. If the point $q_4=Q$ then we can pick $q_4'$ to be the fourth point of $\langle q_3, q_4,q_5\rangle - p^\perp$ and decompose the pentagon into the quadrangle $q_4'^-, q_5^+,q_4^-,q_3^+$ and the pentagon $q_1^-$, $q_2^-$, $q_3^+$, $q_4'^-$, $q_5^+$. We therefore can assume that $q_4\ne Q$.
Moreover modifying this pentagon, if necessary, by the quadrangle
 $q_2^-,q_3^+,q_4^-,Q^+$, we may assume that $q_3=Q$.
But then $q_3^+$ is collinear to $q_1^-$ as well and we can decompose the pentagon into the triangle $q_1^-,q_2^-,q_3^+$ and the quadrangle
 $q_1^-,q_3^+,q_4^-,q_5^+$.
\qed

\bco If $|\FF|=2$ and $n=6$ then the fundamental group of $\Pi(p,H)$ is $\ZZ/2\ZZ$.
\eco
\pf
The geometry $\Pi(p,H)$ has a $2$-cover $\spi$ by Lemma~\ref{lem:2-cover}.
This $2$-cover is simply connected by Lemma~\ref{lem:spi simply connected}.
Therefore $\spi$ is the universal cover of $\Pi(p,H)$ and the fundamental group
 is $\ZZ/2\ZZ$.
\qed
\section{The group action}\label{section:group action}
\subsection{The classical amalgam}\label{sec:classical amalgam}
In this section we shall prove that $G=\Sp(V)$ is the universal completion of subgroups of small rank.
Our first aim is to prove Theorem~\ref{thm:maxparam} using Tits' Lemma~\ref{lem:tits lemma}. To this end, we proved in Section~\ref{section:simple connectedness} that the geometry $\Gamma$ is simply connected. The other result we shall need is the following.

\ble\label{lem:flag-transitive}
The symplectic group $\Sp(V)$ acts flag-transitively on $\Gamma$.
\ele
\pf
Let $F_1$ and $F_2$ be two flags of the same type.
We prove the lemma by induction on $|F_1|=|F_2|$.
Let $M_i$ be the object of maximal type in $F_i$ for $i=1,2$.
These objects have the same isometry type since their dimensions are equal.
Since neither of them intersects $R=\Rad(V)$, by Witt's theorem (See e.g.~\cite{Ta1992})  there is an isometry $g\in \Sp(V)$ with $gM_1=M_2$.

By induction, there is an element in $\Sp(M_2)$ sending $gF_1$ to $F_2$.
Again by Witt's theorem, this local isometry can be extended to an isometry $h\in\Sp(V)$. 
Thus, the element $hg\in\Sp(V)$ sends $F_1$ to $F_2$, as desired.
\qed

\medskip
\Definition

Our next aim is to describe the stabilizers of flags of $\Gamma$ in $G=\Sp(V)$.
Let $I=\{1,2,\ldots,n-1\}$.
Let $C=\{C_i\}_{i\in I}$ be a chamber of $\Gamma$.
For any $J\sbe I$, let $R_J$ be the $J$-residue on $C$ and let $F_J$ be the flag of cotype $J$ on $C$.
We set
$$\begin{array}{ll}
P_J&=\Stab_G(R_J)\\
 B&=\Stab_G(C).\\
 \end{array}$$
Note that $P_J=\Stab_G(F_J)$ and $B=P_\emptyset = \bigcap_{J\sbe I}P_J$.
The group $B$ is called the {\dfn Borel group} for the action of $G$ on $\Gamma$.

\bth\label{thm:maxparam}
The group $\Sp(V)$ is the universal completion of the amalgam of the maximal parabolic subgroups for the action on $\Gamma$.
\eth
\pf
The group $G=\Sp(V)$ acts flag-transitively on $\Gamma$ by Lemma~\ref{lem:flag-transitive}.
By Proposition~\ref{prop:gamma simply connected}, $\Gamma$ is simply connected and now the result follows from
 Tits' Lemma~\ref{lem:tits lemma}.
\qed

\bco\label{cor:parabolics flag-transitive}
The parabolic subgroup $P_J$ acts flag-transitively on the residue $R_J$.
\eco
\pf
This follows immediately from Lemma~\ref{lem:flag-transitive} and the definition of $P_J$ and $R_J$.
\qed

\bth\label{thm:maxparminparam}
Assume $|\FF|\ge 3$. Let $J\sbe I$ with $|J|\ge 3$.
Then, the parabolic $P_J$ of $\Sp(V)$ is the universal completion of the amalgam  $\{P_{J-\{j\}}\mid j\in J\}$
 of rank $(r-1)$ parabolic subgroups contained in $P_J$.
\eth
\pf
The group $P_J$ acts flag-transitively on the residue $R_J$ by Corollary~\ref{cor:parabolics flag-transitive}.
By Proposition~\ref{prop:residues simply connected} the residue $R_J$ is simply connected if $|J|\ge 3$.
Again the result follows from Tits' Lemma~\ref{lem:tits lemma}.
\qed

\bco\label{cor:minparam}
Let $|\FF|\ge 3$ and $\dim(V)\ge 4$. Then, $\Sp(V)$ is the universal completion of the amalgam $\cA_{\le 2}=\{P_J\mid |J|\le2\}$ of rank $\le 2$ parabolic subgroups for the action on $\Gamma$.
\eco
\pf
This follows from Theorems~\ref{thm:maxparam}~and~\ref{thm:maxparminparam} by induction on the rank.
\qed

\medskip

The remainder of this paper is devoted to replacing even the rank $\le 2$ parabolics in the amalgamation result above
 by smaller groups.

\subsection{Parabolic subgroups}\label{subsection:parabolics}
In this section we analyze the parabolic subgroups of rank $\le 2$ of $\Sp(V)$, where $V$ is non-degenerate over a field $\FF$ with $|\FF|\ge 3$. These are the groups in the amalgam $\cA_{\le 2}$ of Corollary~\ref{cor:minparam}.
In  order to study these groups in some detail we will use the following setup.
Let  $n=2r$, and let $\cH=\{e_i,f_j\mid 1\le i\le r, 1\le j\le r\}$ be a hyperbolic basis corresponding to $C$
 as in Corollary~\ref{cor:hyperbolic bases and chambers}.
That is, if we relabel $\cH$ such that
$$\begin{array}{lll}
h_{2i-1}&=e_i & \mbox{ for }1\le i\le r,\\
h_{2j}  &=f_j & \mbox{ for }1\le j\le r,\\
\end{array}$$
then, $C=\{C_l\}_{l\in I}$, where $C_l=\langle h_1,h_2,\ldots,h_l\rangle_V$.
For $i=1,2,\ldots,r$, let $H_i=\langle e_i,f_i\rangle_V$
We will use this setup throughout the remainder of the paper.

Let
$$E=\left(\begin{array}{@{}rr@{}} 0 & 1 \\ -1 & 0 \\\end{array}\right).$$
The matrix defining the symplectic form with respect to the basis $\cH$ is
$$S=
                   \left(\begin{array}{@{}cccc@{}} E       & 0      & \cdots & 0      \\
                                                 0       & \ddots & \ddots & \vdots \\
                                                 \vdots  & \ddots & \ddots & 0      \\
                                                 0       & \cdots & 0      & E      \\
                                                 \end{array}\right)
                   $$

Since the Borel group is contained in all parabolic subgroups, even the ones of rank $\le 2$, we must know exactly what it looks like.

\ble\label{lem:borel}
We have
$$B\cong (\FF\rtimes\GL_1(\FF))^{(n/2)}$$
where, for $j=1,2,\ldots,r$, $\FF\rtimes \GL_1(\FF)$ is realized on $\langle e_j,f_j\rangle$ as
$$B_j=\left\{\left(\begin{array}{@{}ll@{}} a_j & b_j \\ 0 & a_j^{-1}\end{array}\right)\mid a_j\in\FF^*,b_j\in \FF\right\}.$$
Also, the kernel of the action of $G$ on $\Gamma$ is
 $H=\{\pm 1\}$
\ele

\pf
First note that if $n=2$, then $G=\SL_2(\FF)$ and the stabilizer of $C$ is the usual Borel group
$$B_1=\left\{\left(\begin{array}{@{}ll@{}} a_1 & b_1 \\ 0 & a_1^{-1}\end{array}\right)\mid a_1\in\FF^*,b_1\in \FF\right\}.$$
Clearly, $B_1\cong \FF\rtimes \GL_1(\FF)$.

Now let $n\ge 3$.
We observe that, if $g\in G$ stabilizes $C$, then it also stabilizes $C_l^\perp\cap C_m$ for any $1\le l\le m\le n-1$.
Hence, $g$ stabilizes the subspaces spanned by $e_i$ and it stabilizes the subspaces spanned by
 $\{e_i,f_i\}$ for all $1\le i\le r$.
Thus

$$g=               \left(\begin{array}{@{}cccc@{}} g_1       & 0      & \cdots & 0      \\
                                                 0       & \ddots & \ddots & \vdots \\
                                                 \vdots  & \ddots & \ddots & 0      \\
                                                 0       & \cdots & 0      & g_r      \\
                                                 \end{array}\right)
$$
where if  $j=1,2,\ldots,r$ we have
$$g_j=\left(\begin{array}{@{}ll@{}} a_j & b_j \\ 0 & a_j^{-1}\end{array}\right)$$
for some $a_j\in\FF^*$ and $b_j\in \FF$.

The kernel of the action is given by $\bigcap_{x\in G} xBx^{-1}$.
This means that if $g$ is described as above, then, for all $1\le j\le r$ we have
$b_j=0$  and $a_j=a_j^{-1}=a$ for some fixed $a\in \FF^*$. This means $a=\pm 1$.  The result follows.
\qed

\medskip

\medskip

We now focus on the parabolic subgroups of rank $\le 2$.
Since we want to make a distinction between the various rank $\le 2$ parabolic subgroups, we shall give them individual names.

\Definition\label{dfn:sPij sSij sQij}
We assign the following names to the various rank $\le 2$ parabolic subgroups:
$$\begin{array}{lll}
S_j &= P_{2j-1} & \mbox { for } 1\le j\le r\\
M_i &=P_{2i}    & \mbox{ for } 1\le i \le r-1 \\
S_{i\, j}&=P_{2i-1,2j-1} & \mbox{ for } 1\le i<j\le r\\
M_{i\, j}&=P_{2i,2j}       & \mbox{ for }1\le i<j\le r-1\\
Q_{i\, j}&= P_{2i,2j-1}   & \mbox{ for }1\le i\le r-1\mbox{ and } 1\le j\le r\\
\end{array}$$
Thus the collection of groups in the amalgam $\cA_{\le 2}$ is
 $\{M_i,S_j, S_{j\, l},M_{i\, k},Q_{i\, j}\mid 1\le i,k\le r-1\mbox{ and } 1\le j,l\le r\}$.

\Definition\label{dfn: abstract S&M}
In order to describe the groups in $\cA_{\le 2}$ abstractly and as matrix groups we
 define the following matrix groups:

$$
\begin{array}{@{}ll}
M&=\left\{\left.\left(\begin{array}{@{}cccc@{}} a_1 & b_1         & 0 & w \\
                                                               0      & a_1^{-1} & 0 &  0      \\
                                                               0      & a_1^{-1} w a_2 & a_2 & b_2\\
                                                               0      & 0                               & 0     & a_2^{-1}\\
                                                               \end{array}\right)\right| a_1,a_2\in\FF^*,b_1,b_2,w\in \FF\right\},
\\
S&=\SL_2(\FF)\cong \Sp_2(\FF)=\{\left(\begin{array}{@{}cc@{}} a & b\\
                                                               c      & d\\
                                                               \end{array}\right)\mid a,b,c,d\in\FF,\mbox{ where }ad-bc=1\},
\\
M_*&=
\left.\left\{
\left(\begin{array}{@{}cc|cc|cc@{}} a_1 & b_1 & 0 & w_1 & 0 & 0 \\
                                  0 & a_1^{-1} & 0 & 0 & 0 & 0 \\
                                  \hline
                                  0 & a_1^{-1}w_1a_2 & a_2 & b_2 & 0 & w_3 \\
                                   0  & 0   & 0 & a_2^{-1} & 0 & 0  \\
                                  \hline
                                  0 & 0 & 0  & a_2^{-1}w_3a_3    & a_3       & b_3  \\
                                   0 & 0 & 0   & 0 & 0           &  a_3^{-1} \\
      \end{array}
\right)
\right|
\begin{array}{l}
   a_1,a_2,a_3\in \FF^* \\
   b_1,b_2, b_3,w_1,w_3\in\FF\\
\end{array}
\right\},
\\
Q_-&=
\left.\left\{
\left(\begin{array}{@{}cccc@{}} a_1 & b_1 & 0 & w_1 \\
                                  c_1 & d_1 & 0 & w_2\\
                                  v_2 & v_1 & a_2 & b_2\\
                                   0  & 0   & 0 & a_2^{-1}  \\
      \end{array}
\right)
\right|
\begin{array}{l@{}l}
   a_1,b_1,c_1,d_1,w_1,w_2,v_1,& v_2, a_2,b_2\in\FF\mbox{ such that}\\
         a_1d_1-b_1c_1&=1\\
         a_1w_2-c_1w_1+v_2a_2^{-1}&=0\\
         b_1w_2-d_1w_1+v_1a_2^{-1}&=0\\
\end{array}
\right\},
\\
\mbox{and}&\\
Q_+&=
\left.\left\{
\left(\begin{array}{@{}cccc@{}}
                                  a_2     & b_2 & v_2    & v_1 \\
                                  0     & a_2^{-1}     & 0        & 0      \\
                                  0     & w_1 & a_1     & b_1\\
                                   0    & w_2 & c_1     & d_1  \\
      \end{array}
\right)
\right|
\begin{array}{l@{}l}
   a_1,b_1,c_1,d_1,w_1,w_2,v_1,& v_2, a_2, b_2\in\FF\mbox{ such that}\\
         a_1d_1-b_1c_1&=1\\
         a_1w_2-c_1w_1+v_2a_2^{-1}&=0\\
         b_1w_2-d_1w_1+v_1a_2^{-1}&=0\\
\end{array}
\right\}.
\end{array}
$$

\ble\label{lem:structure parabolics}
\begin{itemize}
\AI{a}
For all indices that apply, we have the following isomorphisms:
$$\begin{array}{ll}
S_j &\cong S\times \Pi_{i\ne j} B_i, \\
M_i &\cong M\times\Pi_{j\ne i,i+1} B_j, \\
S_{i\, j}&=\langle S_i, S_j\rangle_G, \\
M_{i\, j}&=\langle M_i, M_j\rangle_G, \\
Q_{i\, j}&=\langle M_i, S_j\rangle_G.\\
\end{array}$$
Moreover, we have the following isomorphisms,
\AI{b}
$$M_{i\, j}\cong\left\{
\begin{array}{ll}
M_*\times \Pi_{k\ne i,i+1,j,j+1} B_k&\mbox{ if }|i-j|=1\\
M\times M \times\Pi_{k\ne i,i+1,j,j+1} B_k &\mbox{ if } |i-j|\ge 2\\
\end{array}\right.,$$
\AI{c}
$$S_{i\, j}\cong (S\times S)\times  \Pi_{k\ne i,j} B_k,$$
\AI{d}
If $j\not\in\{i,i+1\}$, then
$$Q_{i\, j}\cong (M\times S)\times \Pi_{k\ne i,i+1,j} B_k.$$
Furthermore,
$$
Q_{i\, i}\cong Q_-\times\Pi_{k\ne i,i+1} B_k,
$$
and
$$
Q_{i\, i+1}\cong
Q_+\times\Pi_{k\ne i,i+1} B_k.
$$
\end{itemize}
\ele
\pf
(a) By definition $S_i$ is the stabilizer of the flag of type $I-\{2i-1\}$ on the standard chamber $C$.
Therefore it acts on $H_j$ as $B_j$ for all $j\ne i$ and it stabilizes and acts as $S\cong\Sp_2(\FF)$ on $H_i$.
Similarly, by definition, $M_i$ acts on $H_j$ as $B_j$ for all $j\ne i,i+1$ and it stabilizes $\langle e_i\rangle$ and
 $\langle e_i,f_i,e_{i+1}\rangle$ and $H_i\oplus H_{i+1}$. Therefore it acts as $M$ on $H_i\oplus H_{i+1}$.
The last three equalities follow from the fact that rank $2$ residues are connected.

The remaining isomorphisms follow from similar considerations.
The parabolic subgroup under consideration acts as $B_j$ on $H_j$ for all but $2$,$3$, or $4$ values of $j$.
Then the action on the subspace generated by the remaining $H_j$ determines the non-Borel part ($S$, $M$, $M_*$, $Q_-$, $Q_+$) of the group.
\qed

\section{The slim amalgam}\label{section:slim}
In this section we define a slim version $\pcA$ of the amalgam $\cA_{\le 2}$ by eliminating
 a large part of the Borel group from each of its groups.
More precisely, the collection of groups in $\pcA$ is $\{\pX\mid X\in \cA_{\le 2}\}$, where
 $\pX$ is given below in Definitions~\ref{dfn:sPi sSj} and~\ref{dfn:pA}.
Note that the groups $\pX$ are subgroups of $G$ and the inclusion maps for the amalgam $\pcA$ are the inclusion maps induced by $G$.
In Section~\ref{section:concrete amalgam} we shall construct an abstract version of $\pcA$, whose groups are not considered to be subgroups of $G$.

\Definition\label{dfn:sPi sSj}
For $1\le i\le r-1$, the group $\pM_i$ fixes every vector in $H_k$, for every $k\ne i,i+1$ and
a generic element acts on $H_i\oplus H_{i+1}$ as
$$m(b_1,w,b_2)=\left(\begin{array}{@{}cccc@{}} 1 & b_1 & 0 & w \\
                                      0 & 1 & 0 & 0 \\
                                      0 & w & 1 & b_2 \\
                                      0 & 0 & 0 & 1 \\
            \end{array}\right), \mbox{ where } w,b_1,b_2\in \FF.
$$
For $1\le j\le r$, the group $\pS_j$ fixes every vector in $H_k$, for every $k\ne j$ and acts on $H_j$ as $\Sp(H_j)\cong\Sp_2(\FF)$.
A generic element of $\pS_j$ is denoted
$$s(a,b,c,d)= \left(\begin{array}{@{}ll@{}} a & b \\
                                      c & d\\
            \end{array}\right),$$
where the matrix defines the action on $H_j$ with respect to the basis $\{e_j,f_j\}$.

\medskip


\Definition\label{dfn:pA}
For $1\le i<j\le r$, let $$\pS_{i\, j}=\langle\pS_i,\pS_j\rangle_G.$$
For $1\le i<j\le r-1$, let $$\pM_{i\, j}=\langle\pM_i,\pM_j\rangle_G.$$
For $1\le i\le r-1$ and $1\le j\le r$, let $$\pQ_{i\, j}=\langle\pM_i,\pS_j\rangle_G.$$

\Definition\label{dfn:pU pB}
For any $i=1,2,\ldots,r$, we set
$$\pU_i=\left\{
\begin{array}{ll}
\pM_i\cap \pS_i &\mbox{ if }1\le i<r\\
\pM_{i-1}\cap \pS_i & \mbox{ if }i=r\\
\end{array}\right.,$$
and
$$\pB_i=N_{\pS_i}(\pU_i).$$
Moreover,
$$\pB=\langle \pB_1,\ldots,\pB_r\rangle\cong\pB_1\times\cdots\times\pB_r.$$

\Definition\label{dfn: abstract slim S&M}
In order to describe the groups in $\pcA$ abstractly and as matrix groups we
 define the following matrix groups:
$$
\begin{array}{@{}ll}
\pM&=\left.\left\{\left(\begin{array}{@{}cccc@{}} 1 & b_1         & 0 & w \\
                                                               0      & 1 & 0 &  0      \\
                                                               0      & w & 1 & b_2\\
                                                               0      & 0         & 0     & 1\\
                                                               \end{array}\right)\right| b_1,b_2,w\in \FF\right\},
\\
\pS&=S,\\
\pM_*&=
\left.\left\{
\left(\begin{array}{@{}cc|cc|cc@{}} 1 & b_1 & 0 & w_1 & 0 & 0 \\
                                  0 & 1 & 0 & 0 & 0 & 0 \\
                                  \hline
                                  0 & w_1 & 1 & b_2 & 0 & w_3 \\
                                   0  & 0   & 0 & 1 & 0 & 0  \\
                                  \hline
                                  0 & 0 & 0  & w_3   & 1       & b_3  \\
                                   0 & 0 & 0   & 0 & 0           &  1 \\
      \end{array}
\right)
\right|
\begin{array}{l}
      b_1,b_2, b_3,w_1,w_3\in\FF\\
\end{array}
\right\},\\
\pQ_-&=
\left.\left\{
\left(\begin{array}{@{}cccc@{}} a_1 & b_1 & 0 & w_1 \\
                                  c_1 & d_1 & 0 & w_2\\
                                  v_2 & v_1 & 1 & b_2\\
                                   0  & 0   & 0 & 1  \\
      \end{array}
\right)
\right|
\begin{array}{l@{}l}
   a_1,b_1,c_1,d_1,w_1,w_2,&v_1, v_2, b_2\in\FF\mbox{ such that}\\
         1 = a_1d_1-b_1c_1&\\
    v_2=    -a_1w_2+c_1w_1&\\
      v_1=   - b_1w_2+d_1w_1&\\
\end{array}
\right\},
\\
\mbox{and} & \\
\pQ_+&=
\left.\left\{
\left(\begin{array}{@{}llll@{}}
                                  1     & b_2 & v_2    & v_1 \\
                                  0     & 1     & 0        & 0      \\
                                  0     & w_1 & a_1     & b_1\\
                                   0    & w_2 & c_1     & d_1  \\
      \end{array}
\right)
\right|
\begin{array}{l@{}l}
   a_1,b_1,c_1,d_1,w_1,w_2,&v_1, v_2, b_2\in\FF\mbox{ such that}\\
         1=a_1d_1-b_1c_1&\\
        v_2=    -a_1w_2+c_1w_1&\\
      v_1=   - b_1w_2+d_1w_1&\\
\end{array}
\right\}.
\end{array}$$

\ble\label{lem:sPi and sS}
For $i=1,2,\ldots,r-1$ and $j=1,2,\ldots,r$, the sets $\pM_i$ and $\pS_j$ are subgroups of $\Sp(V)$ and
\begin{itemize}
\AI{a} $\pM_i\cong \pM\cong \FF^3$,
\AI{b} $\pS_j\cong \pS\cong \Sp_2(\FF)$.
\end{itemize}
\ele
\pf
For $m=m(b_1,w,b_2)$ one verifies easily that $(m e_i, m f_i)$ and $(m e_{i+1}, m f_{i+1})$
 are two orthogonal hyperbolic pairs. Hence $m$ is a symplectic matrix.
Clearly $\pS_j=\Stab_{\Sp(V)}(H_j)$ is a subgroup of $\Sp(V)$.

(a) It is straightforward to check that $\pM_i$ is an abelian group isomorphic to $\FF^3$.

(b) This is true by definition.
\qed

\medskip

\ble\label{lem:sPij sQij sSij}
For all indices $i\ne j$ that apply, we have
\begin{itemize}
\AI{a}
$$\pM_{i\, j}\cong\left\{
\begin{array}{ll}
\pM_* &\mbox{ if }|i-j|=1\\
\pM_i\times\pM_j &\mbox{ if } |i-j|\ge 2\\
\end{array}\right.,$$
\AI{b}
$$\pS_{i\, j}\cong\pS_i\times \pS_j.$$
\AI{c}
If $j\not\in\{i,i+1\}$, then
$$\pQ_{i\, j}\cong\pM_i\times \pS_j.$$
Furthermore,
$$\begin{array}{ll}
\pQ_{i\, i}&\cong\pQ_-,\\
\pQ_{i\, i+1}&\cong\pQ_+.\\
\end{array}
$$
\end{itemize}
\ele
\pf
Part (a) and (b) are straightforward.
As for part (c), if $j\not\in\{i,i+1\}$, then clearly $\pM_i$ and $\pS_j$ commute and intersect trivially.

We now turn to the cases $\pQ_{i\, i}$ and $\pQ_{i\, i+1}$.
First note that conjugation by the permutation matrix that switches $(e_i,f_i)$ and $(e_{i+1},f_{i+1})$, interchanges
 $\pS_i$ and $\pS_{i+1}$, but fixes every element in $\pM_i$.
Thus it suffices to prove the claim for $\pQ_{i\, i}$.

We consider $\pQ_-$ to be represented as a matrix group with respect to the basis $\{e_i,f_i,e_{i+1},f_{i+1}\}$.
We shall now prove that with this identification $\pQ_-=\pQ_{i\, i}$.
To this end we show that $\pQ_-$ is the stabilizer of the vector $e_{i+1}$ in $\Sp(H_i\oplus H_{i+1})$.
It is clear from the shape of the third column that $\pQ_-$ stabilizes $e_{i+1}$.
On the other hand, any matrix $A$ in $\Sp(H_i\oplus H_{i+1})$ stabilizing $e_{i+1}$ must have such a
 third column.
It must also have zeroes in the last row as in $\pQ_-$ since in fixing $e_{i+1}$ it must also stabilize $e_{i+1}^\perp$.
Any such matrix $A$ must satisfy the conditions on the entries as indicated in the description of $\pQ_-$ since
$\{Ae_i,Af_i, e_{i+1},Af_{i+1}\}$ must be isometric to $\{e_i,f_i,e_{i+1},f_{i+1}\}$ in order for $A$ to be symplectic.
Therefore $\pQ_-$ equals this stabilizer and hence is a group.

Clearly $\pQ_-$ contains $\pS_i$ and $\pM_i$. We now show that $\langle \pS_i,\pM_i\rangle=\pQ_-$.
We note that if $m=m(0,w,b_2)$ and $s=s(a,b,c,d), s'=s(a_1,b_1,c_1,d_1)$ then $sms^{-1}s'$ is the following matrix:
$$\left(\begin{array}{@{}cccc@{}} a_1 & b_1 & 0 & aw \\
                                  c_1 & d_1 & 0 & cw \\
                                  v_2 & v_1 & 1 & b_2  \\
                                   0  & 0   & 0 & 1   \\
      \end{array}
\right),$$
where $v_1,v_2$ are as in the definition of $\pQ_-$. Therefore all the elements of $\pQ_-$ can be obtained that way.

\qed

\ble\label{lem:structure sPij sQij sSij}
We have
\begin{itemize}
\AI{a}
$$\pM_{i\, j}\cong\left\{
\begin{array}{ll}
\FF^5 &\mbox{ if }|i-j|=1\\
\FF^6 &\mbox{ if } |i-j|\ge 2\\
\end{array}\right.,$$
\AI{b}
$$\pS_{i\, j}\cong\Sp_2(\FF)\times \Sp_2(\FF),$$
\AI{c}
If $j\not\in\{i,i+1\}$, then
$$\pQ_{i\, j}\cong\FF^3\times \Sp_2(\FF).$$

\AI{d}
Taking the labelling as in Definition~\ref{dfn: abstract slim S&M},
let
$$
\begin{array}{@{}ll}
U&=
\left.\left\{
\left(\begin{array}{@{}cccc@{}} 1 & 0 & 0 & 0 \\
                                0 & 1 & 0 & 0\\
                                0 & 0 & 1 & b_2\\
                                0 & 0 & 0 & 1  \\
      \end{array}
\right)
\right|
   b_2\in\FF
\right\}
,\\
V&=
\left.\left\{
\left(\begin{array}{@{}cccc@{}}    1 & 0 & 0 & w_1\\
                                   0 & 1 & 0 & w_2\\
                                  v_2 & v_1& 1 & b_2 \\
                                   0  & 0  & 0 & 1  \\
      \end{array}
\right)
\right|
\begin{array}{l@{}l}
   w_1,w_2,v_1, v_2, b_2\in\FF & \mbox{ such that}\\
         v_2=-w_2&\\
         v_1=w_1&\\
\end{array}
\right\},\\

\pS_i&=
\left.\left\{
\left(\begin{array}{@{}cccc@{}} a_1 & b_1 & 0 & 0 \\
                                  c_1 & d_1 & 0 & 0\\
                                  0 & 0 & 1 & 0\\
                                   0  & 0   & 0 & 1  \\
      \end{array}
\right)
\right|
\begin{array}{l@{}l}
   a_1,b_1,c_1, d_1\in\FF &\mbox{ such that}\\
         a_1d_1-b_1c_1=1 &\\
\end{array}
\right\}
.
\end{array}$$
Then, we have
\begin{itemize}
\AI{i} $U\cong \FF$ and  $U= Z(\pQ_{i\, i})$,
\AI{ii} $V\cong \FF^3$ if $\Char(\FF)=2$, otherwise $V$ is a non-split extension of $U=Z(V)$ by $\FF^2$,
\AI{iii} the action of $\pS_i$ on $V/U$ by conjugation is the natural action from the left of $\Sp_2(\FF)$ on $\FF^2$.
\AI{iv} $\pQ_{i\, i+1}\cong\pQ_{i\, i}\cong \FF. \FF^2\rtimes \Sp_2(\FF)$,
where the first isomorphism is given by the labelling of the entries and the latter follows from (i)-(iii).
\end{itemize}
\end{itemize}
\ele
\pf
(a) and (b): This is straightforward.
(c): This follows from the corresponding part in Lemma~\ref{lem:sPij sQij sSij}.

(d):
The isomorphism  $U\cong \FF$ as well as parts (ii) and (iii) can be verified by straightforward calculation.
One can verify that if $\Char(\FF)\ne 2$ we have $[V,V]=U=Z(V)$.
The isomorphism $\pQ_{i\, i}\cong\pQ_{i\, i+1}$ is given by conjugation as in the proof of Lemma~\ref{lem:sPij sQij sSij}.
Clearly $\langle V,\pS_i\rangle\cong V\rtimes \pS_i$ with the action described above.
It is easy to see that $\pM_i\le \langle V, \pS_i\rangle$.
Therefore $\pQ_{i\, i}=\langle \pM_i,\pS_i\rangle\le \langle V,\pS_i\rangle\le \pQ_{i\, i}$.

Finally, one verifies directly that $\pS_i$ acts trivially on $U$ and by (iii) acts fixed-point freely on $V/U$.
Since $\pQ_{i\, i}=V\rtimes \pS_i$ it now follows that $Z(\pQ_{i\, i})=U$.
By (ii) we have $Z(V)=U$ if $\Char(\FF)\ne 2$.
%
%
%
\qed

\medskip

\section{The concrete amalgam}\label{section:concrete amalgam}

\Definition We will define a concrete amalgam $\scA$. Its set of subgroups is $\{\sM_{i\, j},\sQ_{i\, k},\sS_{k\, l},\sM_i,\sS_k\mid 1\le i,j\le r-1, 1\le k,l\le r\}$, where for each $X\in \cA_{\le 2}$,  $\sX$ is a copy of $\pX$ and
 the inclusion homomorphisms are as follows:
$$\begin{array}{ll}
\varphi^{P}_{i,\{i,j\}}\colon &\sM_i\to\sM_{i\, j}\\
\varphi^{P}_{j,\{i,j\}}\colon &\sM_j\to\sM_{i\, j}\\
\varphi^{S}_{k,\{k,l\}}\colon &\sS_k\to\sS_{k\, l}\\
\varphi^{S}_{l,\{k,l\}}\colon &\sS_l\to\sS_{k\, l}\\
\varphi^{PQ}_{i,\{i,k\}}\colon &\sM_i\to\sQ_{i\, k}\\
\varphi^{SQ}_{k,\{i,k\}}\colon &\sS_k\to\sQ_{i\, k}.\\
\end{array}$$
These inclusions are given by the presentations of $\pX$ as matrix groups as given in
 Definitions~\ref{dfn:sPi sSj} and~\ref{dfn:pA}.
We denote the universal completion of $\scA$ by $\sG$.

\Definition
We now define the map $\pi\colon\scA\to G$.
For any $X\in\cA_{\le 2}$, it identifies $\sX$ with its isomorphic copy,  $\pX$, in $G$.
Thus the image of $\scA$ under $\pi$ is $\pcA$.

\ble\label{lem:pi = iso}
For any element $X\in\cA_{\le 2}$, the map
 $\pi\colon \sX\to\pX$ is an isomorphism.
\ele
\pf
This is true by the definition of $\pi$.
\qed

\medskip

It then follows that the map $\pi$ extends to a surjective map from the universal cover $G^\circ$ to $G$.

\ble\label{lem:intersection of Pi}
If $i\ge2$ then $\sS_i\cap \sM_i=\sS_i\cap\sM_{i-1}=\sM_i\cap \sM_{i-1}$.
\ele
\pf
Note that the above are true for the images under $\pi$. By definition the group $\pM_i$ stablizes the spaces $H_k$ if $k\ne i,i+1$ and it also stabilizes $e_i$ and $e_{i+1}$. Moreover $\pS_i$ stabilizes all the $H_k$ for $k \ne i$. Therefore  $\pS_i\cap \pM_i$ stabilizes all vectors in $H_k$ for $k\ne i$ and it also stabilizes $e_i$. Similarly for $\pS_i\cap\pM_{i-1}$ and $\pM_i\cap \pM_{i-1}$.
We notice that  $\pi$ is an isomorphism when restricted to $\cA^\circ$ and so the conclusion follows.
\qed

\Definition\label{dfn:sU sB}
For $i=1,2,\ldots,r$, we define the following subgroups of $\sS_i\in\scA$:
$\sU_i$ is the common intersection from Lemma~\ref{lem:intersection of Pi} and
$\sB_i=N_{\sS_i}(\sU_i)$.
Furthermore, we set
 $\sB=\langle\sB_i\mid i=1,2,\ldots,r\rangle$ as a subgroup of $G^\circ$.

\ble
For $i=1,2,\ldots, r$, we have
\begin{itemize}
\AI{a} $\pU_i=\pi(\sU_i)$ and $\pB_i=\pi(\sB_i)$.
\AI{b} $$\sU_i=\left\{
\left(\begin{array}{@{}ll@{}} 1 & b_1 \\
                              0 & 1 \\
      \end{array}
\right)
\Big|
\begin{array}{l@{}l}
   b_1\in\FF\\
\end{array}
\right\}
\cong\FF.$$

\AI{c} $$\sB_i=\left\{
\left(\begin{array}{@{}ll@{}} a_1 & b_1 \\
                                  0 & d_1 \\
      \end{array}
\right)
\Big|
\begin{array}{l@{}l}
   a_1,b_1,d_1\in\FF\mbox{ such that }\\
   a_1d_1=1\\
\end{array}
\right\}
\cong\FF\rtimes \FF^*.$$

\end{itemize}

\ele
\pf
By Definition~\ref{dfn:pU pB} and Lemma~\ref{lem:intersection of Pi} and the fact that
 $\pi$ is an isomorphism between $X^\circ$ and $X^\pi$ for every $X\in\cA_{\le 2}$ we have
  $\pi(\sU_i)=\pU_i$. The equality $\pi(\sB_i)=\pB_i$ follows directly from Definitions~\ref{dfn:pU pB}~and~\ref{dfn:sU sB}.
Parts (b) and (c) are verified readily.
\qed

\ble \label{lem:B normalizes}
For $X\in\cA_{\le 2}$, $\sB$ normalizes $\sX$.
The action of $\sB_i$ on the rank $2$ parabolics $\sQ_{i\, i}$, $\sQ_{j\, i}$ ($j\ne i$)
 as well as the action on $\sS_{i\, j}$ is given by conjugation.
Therefore,
\begin{itemize}
\AI{a} $[\sB_i, \sS_j]=1$ if $i\ne j$ and $\sB_i$ acts on $\sS_i$ as inner automorphisms.
\AI{b} $[\sB_i,\sM_j]=1$ if $j \ne i,i-1$. Moreover, $\sB_i$ acts on $\sM_i$ as the conjugation in $\sQ_{i\, i}$ and $\sB_{i+1}$ acts on $\sM_i$ as the conjugation in $\sQ_{i\, i+1}$.
\AI{c} The action of $\sB_i$ on the rank $2$ parabolics $\sM_{i\, j}$ is given by the action on
 its subgroups $\sM_i$ and $\sM_j$.
\end{itemize}
\ele
\pf
Note that the actions are as above for the $\pB_i\le \pS_i$ acting on the various $\pX$. Moreover the map $\pi$ is an isomorphism when restricted to the various $\sX\in\cA^\circ$.
Therefore, the action of $\sB_i$ on any (subgroup of) a rank $2$-parabolic $\sX$ containing $\sS_i$ is the same
 as the action of $\pB_i$ on the (corresponding subgroup of) the corresponding rank $2$-parabolic $\pX$.
This explains the action of $\sB_i$ on the following groups:
(1) $\sS_{i\, j}$ and its subgroups $\sS_i$, $\sS_j$,
(2) $\sQ_{i\, i}$ and its subgroup $\sM_i$,
(3) $\sQ_{j\, i}$ and its subgroup $\sM_j$ ($j\ne i$).
This settles parts (a) and (b).

Part (c) follows immediately since the groups in $\cA^\circ$ are all embedded in $\sG$.
\qed

\ble\label{lem:structure sB}
\begin{itemize}
\AI{a}
The group $\sB$ is the internal direct product
 $$\sB=\sB_1\times \sB_2\times\cdots\times \sB_r,$$
\AI{b}
 $$\pi\colon\sB\to\pB\mbox{ is an isomorphism.}$$
\end{itemize}
\ele
\pf
For any $1\le i<j\le r$, we have $\sB_i\le \sS_i$ and since $\sS_{i\, j}=\sS_i\times \sS_j$ we have
 $\langle \sB_i,\sB_j\rangle=\sB_i\times \sB_j$.
\qed

\ble\label{lem: pBpX=X}
For $X\in\cA_{\le 2}$ we have
 $\langle B^\pi,X^\pi\rangle_G=X$.
\ele
\pf
This is an easy calculation inside $G$.
\qed

\ble\label{lem:pi preserves intersection}
For any $X\in\cA_{\le 2}$, $\pi(\sX\cap\sB)=X^\pi\cap B^\pi$.
\ele
\pf
The inclusion $\sbe$ is trivial.
We now prove $\spe$.
We do this case by case for any $X\in\cA$.

Let $X=S_i$.
Then $\pX\cap\pB=\pB_i$ and since $\sB_i\le \sS_i\cap\sB$, we find $\pi(\sX\cap\sB)\spe \pi(\sB_i)=\pB_i$.

Let $X=P_i$.
Then $\pX\cap\pB=\pU_i\times \pU_{i+1}$ and since $\sU_i\times \sU_{i+1}\le \sM_i\cap\sB$, we find
$\pi(\sX\cap\sB)\spe \pi(\sU_i\times\sU_{i+1})=\pU_i\times \pU_{i+1}$.

Let $X=S_{i\, j}$.
Then $\pX\cap\pB=\pB_i\times \pB_j$ and since $\sB_i\times \sB_j\le \sS_{i\, j}\cap\sB$, we find
$\pi(\sX\cap\sB)\spe \pi(\sB_i\times\sB_j)=\pB_i\times \pB_j$.

Let $X=Q_{i\, j}$.
Then $\pX\cap\pB=\langle\pU_i, \pU_{i+1}, \pB_j\rangle$ and since $\langle\sU_i, \sU_{i+1}, \sB_j\rangle\le \sQ_{i j}\cap\sB$, we find
$\pi(\sX\cap\sB)\spe \pi(\langle\sU_i, \sU_{i+1}, \sB_j\rangle)=\langle\pU_i, \pU_{i+1}, \pB_j\rangle$.

Let $X=P_{i\, j}$.
Then $\pX\cap\pB=\langle\pU_i, \pU_{i+1}, \pU_j, \pU_{j+1}\rangle$ and since $\langle\sU_i,$ $\sU_{i+1}$, $\sU_j$, $\sU_{j+1}\rangle\le \sM_{i\, j}\cap\sB$, we find
$\pi(\sX\cap\sB)\spe \pi(\langle\sU_i, \sU_{i+1}, \sU_j, \sU_{j+1}\rangle)=\langle\pU_i$, $\pU_{i+1}$, $\pU_j$, $\pU_{j+1}\rangle$.
\qed

\mn
Our next aim is to extend the map $\pi^{-1}\colon \pX\to \sX$ for every $X\in\cA_{\le 2}$. To this end, for each such $X$, and all $a\in X$, define $\chi(a)=\pi^{-1}(a)$.

\ble\label{lem:well defined on intersections}
$\chi$ is well-defined on $X^\pi\cap B^\pi$ for all $X\in\cA_{\le 2}$.
\ele

\pf
This follows from Lemma~\ref{lem:pi preserves intersection}
\qed

\mn
Define $\chi$ on $X=B^\pi X^\pi$ for any $X\in\cA_{\le 2}$ as follows:
$\chi(bx)=\chi(b)\chi(x)$.

\ble\label{lem:chi injective}
$\chi$ is well-defined and injective on $X$.
\ele
\pf
Note that if $b,b' \in \pB$ and $x, x'\in \pX$ then $bx=b'x' $ implies $b^{-1}b'=xx'^{-1} \in \pB\cap \pX$. Moreover by Lemma \ref{lem:well defined on intersections} $\chi(b^{-1}b')=\chi(xx'^{-1})$ and so, using that $\chi$ is an isomorphism when restricted to $\pB$ and $\pX$, it follows that $\chi(bx)=\chi(b'x')$ and $\chi$ is well defined.

Also if $bx\in X$ with $\chi(bx)=1$ it follows that $\chi(b)=\chi(x)^{-1} \in \sB\cap\sX$ (because $\chi$ is the inverse of $\pi$ when restricted to $\pB$ and $\pX$). Therefore $b=\pi(\chi(b))$ and $x=\pi(\chi(x))$ are both in $\pi(\sX\cap\sB)=X^\pi\cap B^\pi$. But $\chi$ is a bijection when restricted to $\pX\cap \pB$ and so $b=x^{-1}$ and $\chi$ is injective.
\qed

\ble\label{lem:chi amalgam hom}
$\chi$ is an embedding of the amalgam $\cA_{\le 2}$ into $\sG$.
\ele
\pf
Let $X \in \cA$ then by Lemma~\ref{lem: pBpX=X} and Lemma~\ref{lem:B normalizes} we know that $X=\pB\pX$ and so if $bx, b'x' \in X$ then $\chi(b'x'bx)=\chi(b'b)\chi(b^{-1}x'bx)=\chi(b'b)\chi(b^{-1}x'b)\chi(x)$. Also $\chi(b'x')\chi(bx)=\chi(b')\chi(x')\chi(b)\chi(x)=\chi(b')\chi(b)\chi(b^{-1})\chi(x')\chi(b)\chi(x)$. So we only need to prove that $\chi(b^{-1})\chi(x')\chi(b)=\chi(b^{-1}x'b)$ which is equivalent to the fact that the action of $\pB$ on $\pX$ is the same as the action of $\sB$ on $\sX$. This follows from Lemma~\ref{lem:B normalizes}.
\qed

As a consequence, we find the following result.

\bpr\label{prop:chi surjective}
The map $\chi$ extends to a surjective homomorphism $G\to \sG$ which we also denote by $\chi$.
\epr
\pf
By Lemma~\ref{lem:chi amalgam hom} the map $\chi$ extends to a homomorphism $G\to\sG$ whose image contains
 the subgroup of $\sG$ generated by the subgroups in $\cA^\circ$.
Since $\langle\scA\rangle=\sG$, the conclusion follows.
\qed

\medskip
\paragraph{Proof of Theorem~\ref{thm:main theorem1}.}
By Corollary~\ref{cor:minparam} $G$ is the universal completion of the amalgam $\cA_{\le 2}$.

The map $\pi\colon \cA^\circ\to G$ given by
 $\pi\colon X^\circ\to X^\pi$ for all $X\in\cA_{\le 2}$ extends to surjection
   $\pi\colon G^\circ\to H^\pi$ where $H^\pi$ is the subgroup of $G$ generated by the subgroups in $\cA^\pi$.
By Lemma~\ref{lem: pBpX=X} $H^\pi=G$, and so $\pi\colon G^\circ\to G$ is surjective.

By Proposition~\ref{prop:chi surjective} there is a surjective map $\chi\colon G\to \sG$, which by Lemma~\ref{lem:chi injective} is injective on the subgroups of
 $\cA_{\le 2}$.

The composition $\chi\after \pi\colon \sG\to\sG$ is a surjective homomorphism which is the identity on every subgroup in
 the amalgam $\cA^\circ$.
Since $\sG$ is the universal completion of $\cA^\circ$, the only such map is the identity.

Hence $\pi\colon\sG\to G$ is an isomorphism.
\qed

\bibliographystyle{plain}




\end{document}